\documentclass[12pt]{article}
\usepackage{amssymb}
\usepackage{amsmath}
\usepackage{latexsym}

\newtheorem{theorem}{Theorem}[section]
\newtheorem{proposition}[theorem]{Proposition}

\newtheorem{lemma}[theorem]{Lemma}

\newtheorem{remark}[theorem]{Remark}

\newcommand{\half}{\frac{1}{2}}
\newcommand{\tpi}{2\pi i}
\begin{document}
%%%%%%%%%%%%%%%%%%%%%%%%%%%%%%%%%%%%%%%%%%%%%%%%%%%%%%%%%%%%%
\title{The Szeg\"o Kernel on a Sewn Riemann Surface}
\author{
 Michael P. Tuite and Alexander Zuevsky
\thanks{%
Supported by a Science Foundation Ireland Research Frontiers Programme Grant}
 \\
School of Mathematics,  
Statistics and Applied Mathematics,  
\\
%EndAName
National University of Ireland, 
Galway, Ireland\\
}
\maketitle
%%%%%%%%%%%%%%%%%%%%%%%%%%%%%%%%%%%%%%%%%%%%%%%%%%%%%%%%%%
\begin{abstract}
We describe the Szeg\"o kernel on a higher genus Riemann surface in terms of Szeg\"o kernel data coming from lower genus surfaces via two explicit sewing procedures where either two Riemann surfaces are sewn together or a handle is sewn to a Riemann surface.  
We consider in detail the examples of the Szeg\"o kernel on a genus two Riemann surface formed by either sewing together two punctured tori or by sewing a twice-punctured torus to itself. We also consider the modular properties of the Szeg\"o kernel in these cases.\end{abstract}

%%
%%%%%%%%%%%%%%%%%%%%%%%%%%%%%%%%%%%%%%%%%%%%%%%%%%%%%%%%%%%%%
\section{Introduction}
\label{Sect_Intro}
%%%%%%%%%%%%%%%%%%%%%%%%%%%%%%%%%%%%%%%%%%%%%%%%%%%%%%%%%%%%%%%%%
%%
The purpose of this paper is to provide an explicit description of
the Szeg\"o kernel \cite{Sz},\cite{HS},\cite{Sc},\cite{F1} on a higher genus
Riemann surface in terms of Szeg\"o kernel data coming from lower genus surfaces. We exploit two explicit sewing procedures where either 
two lower genus Riemann surfaces are sewn together or else a handle is sewn to a lower genus Riemann surface. We also consider in some detail the
construction and modular properties of the Szeg\"o kernel on a genus two Riemann surface formed either by sewing two tori together or 
by sewing a handle on to a torus.
This paper is a further development of the theory of partition and 
$n$-point correlation functions on Riemann surfaces for vertex operator algebras (e.g. \cite{FLM},\cite{Ka}) as described 
in \cite{T},\cite{MT1},\cite{MT1},\cite{MT2},\cite{MT3},\cite{MT4},\cite{MTZ}. Our main motivation is to lay the foundations for the explicit construction
of the partition and $n$-point correlation functions for a fermionic vertex operator super algebra on higher genus 
Riemann surfaces \cite{TZ1},\cite{TZ2}. (The central role played by the Szeg\"o kernel for such systems has been long appreciated in theoretical physics \cite{RS}, \cite{R}, \cite{DVFHLS}, \cite{DVPFHLS}). Our ultimate aim is to develop a fully rigorous theory of higher genus partition and $n$-point functions along the lines of Zhu's theory at genus one \cite{Z}. Thus, to name but a few applications, one might eventually study Siegel (sub)group modular invariance, the Freidan-Shenker conjecture concerning the reconstructability of a conformal field theory from the genus $g$ partition function at all genera \cite{FS}, higher genus bosonization etc.
However, this present paper may also be of interest to readers outside the vertex operator algebra community.

%%%%%%%%%%%%%%%%%%%%%%%%%%%%%%%%%%%%%%%%%%%%%%%%%%%%%%%%%%%%%%%%%%%%%%%%%%%%%%%%
\medskip

%%%%%%%%%%%%%%%%%%%%%%%%%%%%%%%%%%%%%%%%%%%%%%%%%%%%%%%%%%%
We begin in Section \ref{Szegokernel} with a review of some basic aspects of the theory of Riemann surfaces 
\cite{FK},\cite{Sp},\cite{F1},\cite{F2},\cite{Mu}. We then define and discuss properties of the Szeg\"o kernel, which is a meromorphic $(\half,\half)$ 
differential with a simple pole structure and prescribed multiplicities on the cycles of the Riemann surface \cite{Sz},\cite{HS},\cite{Sc},\cite{F1}.

\medskip
%%%%%%%%%%%%%%%%%%%%%%%%%%%%%%%
 In Section \ref{Sect_Epsilon_g} we describe the Szeg\"o kernel on a genus $g_{1} +g_2$ Riemann surface $\Sigma^{(g_{1} + g_2)}$ 
obtained by sewing two lower genus Riemann surfaces $\Sigma^{(g_{1} )}$ and  $\Sigma^{(g_2)}$. This is similar to the 
approach of refs.~\cite{Y} and \cite{MT1}, for computing the period matrix and other related structures on $\Sigma^{(g_{1} + g_2)}$ 
in terms of lower genus data. 
Following \cite{MT1}, we refer to this sewing scheme as the \emph{$\epsilon$-formalism} where  
$\epsilon$ is a complex sewing parameter which forms part of
the data according to which the sewing is performed (see Figure~1 below). 
In particular, we introduce an infinite block matrix 
\begin{equation}
\label{qmatrix}
Q=\left( 
\begin{array}{cc}
0 & \xi F_{1}
%%\left( { \theta}_{1}   \atop { \phi }_{1}  \right) 
\\ 
-\xi F_{2}
%%\left( { \theta}_2  \atop { \phi }_2 \right) 
& 0
\end{array}
\right), 
\end{equation}
where $F_{1}$ and $F_2$ are infinite matrices whose entries are
certain weighted moments of the Szeg\"o kernels on $\Sigma^{(g_1)}$ and $\Sigma^{(g_2)}$, respectively, and
$\xi\in \{\pm \sqrt{-1}\}$. The matrix $I-Q$, where $I$ is the infinite identity matrix, plays a crucial role
here (and in the sequel \cite{TZ1}). In particular, we show that $I-Q$ is invertible for small enough
$\epsilon$.  $(I-Q)^{-1}$ then forms part of the expression of the genus $g_1+g_2$ Szeg\"o kernel in 
terms of the lower genus Szeg\"o kernel data as proved in Theorem \ref{Theorem_Shh}. In Theorem \ref{theorem_Det} we further show that the determinant $\det (I -Q)$ is well-defined and is a 
non-vanishing holomorphic function for small enough $\epsilon$. Finally, we describe the 
example of the Szeg\"o kernel on a genus two Riemann surface formed by sewing two tori and 
verify its modular  transformation properties under the modular group which preserves the 
sewing scheme. This example is extensively exploited in \cite{TZ1}. 
\medskip

%%%%%%%%%%%%%%%%%%%%%%%
Section \ref{Sect_Rho_g} is devoted to
development of the corresponding formalism in the case that $\Sigma^{(g+1)}$
of genus $g+1$ is obtained
by self-sewing a handle to a genus $g$ Riemann surface $\Sigma^{(g)}$ with complex sewing parameter $\rho$.
We refer to this as the \emph{$\rho$-formalism}. This case is more technical due to the extra multiplicities on the two new cycles associated with the sewing handle. This leads us to
introduce an analogue  of \eqref{qmatrix}, namely, an infinite matrix 
$T$ whose entries are determined by weighted moments of certain genus $g$ objects related to the Szeg\"o kernel on $\Sigma^{(g)}$ and the new multiplicities.  We show that $I-T$ is invertible for suitably small $\rho$ and in Theorem \ref{Theorem_Shhrho} express the Szeg\"o kernel on  
$\Sigma^{(g+1)}$ in terms of $(I-T)^{-1}$ and other genus $g$ Szeg\"o kernel data. 
In Theorem \ref{theorem_Det_rho} we show that the determinant
 $\det(I-T)$ is well-defined and holomorphic for suitably small $\rho$.
We conclude with two examples of sewing a handle to a Riemann sphere to obtain a torus and sewing a handle to a torus to obtain genus two Riemann surface. The modular transformation properties of the genus two Szeg\"o kernel are also verified under the modular group preserving this $\rho$-sewing scheme. This example will be extensively exploited in \cite{TZ2}.

%%%%%%%%%%%%%%%%%%%%%%%%%%%%%%%%%%%%%%%%%%%%%%%%%%%%%%%%%%%%%%%%%%%%%%
%%%%%%%%%%%%%%%%%%%%%%%%%%%%%%%%%%%%%%%%%%%%%%%%%%%%%%%%%%%%%%%%%%
\section{The Szeg\"o Kernel on a Riemann Surface}  
\label{Szegokernel}
%%%%%%%%%%%%%%%%%%%%%%%%%%%%%%%%%%%%%%%%%%%%%%%%%%%%%%%%%%%%%%%%%%
%%
Consider a compact Riemann surface $\Sigma$ of genus $g$ with canonical
homology cycle basis $a_{1},\ldots, a_{g}, b_{1}, \ldots, b_{g}$. In general there
exists $g$ holomorphic 1-forms $\nu_{i}$, $i=1,\ldots, g$ which we may
normalize  by (e.g. \cite{FK},\cite{Sp}) 
\begin{equation}
\oint_{a_{i}}\nu_{j}=\tpi \delta_{ij}.  \label{norm}
\end{equation}%
The genus $g$ period
matrix $\Omega $ is defined by 
\begin{equation}
\Omega_{ij}=\frac{1}{\tpi }\oint_{b_{i}}\nu_{j}, 
\label{period}
\end{equation}%
for $i,j=1,\ldots, g$. $\Omega$ is symmetric with positive imaginary part i.e. $\Omega\in \mathbb{H}_g$, the Siegel upper half plane. 
The canonical intersection form on cycles  is preserved under the action of the symplectic group  $Sp(2g,\mathbb{Z})$ where
\begin{equation}
\left(\begin{array}{c}
{b} \\ 
{a}
\end{array}\right)
\rightarrow
\left(\begin{array}{c}
\tilde{b} \\ 
\tilde{a}
\end{array}\right)
=
\left(\begin{array}{cc}
A & B \\ 
C & D
\end{array}\right)
\left(\begin{array}{c}
{b} \\ 
{a}
\end{array}\right),\quad \left(\begin{array}{cc}
A & B \\ 
C & D
\end{array}\right)\in Sp(2g,\mathbb{Z}).
\label{eq:modcycle}
\end{equation}
This induces the modular action on $\mathbb{H}_g$ 
\begin{equation}
\Omega\rightarrow \tilde \Omega=\left(A\Omega+B\right)\left(C\Omega+D\right)^{-1}.
\label{eq:modOmega}
\end{equation} 

It is useful to introduce the \emph{normalized differential of the
second kind} defined by \cite{Sp},\cite{Mu},\cite{F1}: 
\begin{equation}
\omega(x,y)\sim  \frac{dxdy}{(x-y)^{2}} \quad \mbox{for } x\sim y, 
\label{omegag}
\end{equation}
for local coordinates $x,y$, with normalization 
$
\int_{a_{i}}\omega(x,\cdot )=0$ for $i=1,\ldots, g$. Using the Riemann bilinear relations, one finds that 
$\nu_{i}(x)=\oint_{b_{i}}\omega(x,\cdot )$.

We also introduce the \emph{normalized
differential of the third kind} 
\begin{equation}
\omega_{p_{2}-p_{1}}(x)=\int_{p_{1}}^{p_{2}}\omega(x,\cdot ),
\label{omp2p1}
\end{equation}
for which $\oint_{a_{i}}\omega_{p_{2}-p_{1}}=0$ and
$\omega_{p_{2}-p_{1}}(x)\sim  \frac{(-1)^{a}}{x-p_{a}}dx$ 
for $x\sim p_{a}$ and $a=1,2$.

We recall the definition of the theta function with real characteristics  e.g. \cite{Mu},\cite{F1},\cite{FK}
\begin{equation}
\label{theta}
\vartheta \left[ {{\alpha} \atop {\beta} }\right] \left( {z}\vert  \Omega  \right)  
=
\sum_{ m \in {\mathbb{Z}}^{g} } 
\exp \left(
%e^{%%
 i \pi (m + {\alpha}).\Omega .(m+{\alpha}) + 
 (m+ {\alpha}). ({{z}+ \tpi  {\beta}}) 
\right),
% },   
%%
\end{equation}
for $\alpha=(\alpha_i),{\beta}=(\beta_i)\in \mathbb{R}^g$, ${z}=(z_i)\in \mathbb{C}^g$ and $i=1,\ldots, g$ with
\begin{eqnarray}
\vartheta \left[ {{\alpha} \atop {\beta} }\right] 
\left(z + \tpi  (\Omega. r+s) \vert  \Omega\right)
&=&
 e^{\tpi  {\alpha}.s} e^{-\tpi  {\beta}.r} e^{-i\pi  r.\Omega.r -r.z}
\vartheta \left[ { {\alpha} \atop {\beta} }\right] \left(z \vert  \Omega\right), \notag\\
\vartheta \left[ {{\alpha+r} \atop {\beta+s} }\right] 
\left(z \vert  \Omega \right) 
&=& 
e^{\tpi  {\alpha}.s}\vartheta \left[ { {\alpha} \atop {\beta} }\right] \left(z \vert  \Omega\right),\label{theta_char_periodicity}  
\end{eqnarray}
for $r,s \in {\mathbb Z}^g$.

There exists a (nonsingular and odd) character $\left[ { \gamma} \atop {\delta} \right]$ such that \cite{Mu},\cite{F1} 
\begin{eqnarray}
\vartheta \left[ { \gamma} \atop {\delta} \right](0\vert  \Omega)=0,\qquad
 \partial _{z_i}\vartheta \left[ { \gamma} \atop {\delta} \right](0\vert  \Omega)\neq 0.
\label{eq:thetadelta}
\end{eqnarray}
Let 
\begin{eqnarray}
\zeta(x) = 
\sum_{i=1}^g 
\partial _{z_i}\vartheta \left[ { \gamma} \atop {\delta} \right] (0\vert  \Omega)\nu_i(x), 
\label{eq:zeta}
\end{eqnarray}
a holomorphic 1-form, and let $\zeta(x)^{\half}$ denote the form of weight ${\half}$ on the double cover $\widetilde\Sigma$ of $\Sigma$. We also refer to $\zeta(x)^{\half}$ as a (double-valued) ${\half}$-form on $\Sigma$. We define the prime form $E (x, y)$ by\footnote{Note that our definition differs from that of refs. \cite{Mu,F1} by a factor of $-1$.}
\begin{equation}
\label{prime_form}
E(x,y) =\frac{ \vartheta \left[ {{ \gamma} \atop {\delta}} \right] 
\left( \int_{y}^{x}\nu \vert \Omega\right) } 
 {\zeta(x)^{\half}\zeta(y)^{\half}}\sim (x-y)dx^{-\half} dy^{-\half} \quad \mbox{for } x\sim y,   
\end{equation}
where $\int_{y}^{x}\nu= (\int_{y}^{x}\nu_i)\in \mathbb{C}^g$. $E(x,y)=-E(y,x)$ is a holomorphic differential form of weight $(-\half,-\half)$ on $\widetilde{\Sigma} \times \widetilde{\Sigma}$.  $E(x,y)$ has multipliers along the $a_i$ and $b_j$ cycles in $x$ given by $1$ and $e^{-i\pi \Omega_{jj}-\int_{y}^{x}\nu_j}$ respectively \cite{F1}. 

The normalized differentials of the
second and third kind can be
expressed in terms of the prime form  \cite{Mu} 
\begin{eqnarray}
\omega(x,y) &=&\partial_{x}\partial_{y}\log E(x,y)dxdy, 
%%
%%%
\label{omprime}
\\
\omega_{p-q}(x) &=&\partial_{x}\log 
\frac{E(x,p)}{E(x,q)}dx.  
\label{omp2p1prime}
\end{eqnarray}
Conversely, we can also express the prime form in terms of $\omega$ by \cite{F2}
\begin{equation}
E(x,y)=\lim_{p\rightarrow x,\ q \rightarrow y}\left[
\sqrt{(x-p)(q-y)}\exp \left (-\half\int_y^x\omega_{p-q} \right )
\right]
dx^{-\half}dy^{-\half}.
\label{eq:Prime_omega}
\end{equation}

We define the Szeg\"o kernel \cite{Sc},\cite{HS},\cite{F1} for $\vartheta \left[ {{\alpha} \atop {\beta}}\right] ( 0\vert  \Omega)\neq 0$ as follows
\begin{equation}
\label{Szegodefn}
 S\left[ {{\theta} \atop {\phi} }\right] (x, y\vert \Omega) =
\frac{ \vartheta \left[ {{\alpha} \atop {\beta}} \right] 
\left( \int_{y}^{x}\nu \vert  \Omega\right) } %%
  {\vartheta \left[ {{\alpha} \atop {\beta}}\right] ( 0\vert  \Omega)  E(x, y)},
\end{equation}
 where $\theta=({\theta}_i),\ \phi=(\phi_i)\in U(1)^n$ for
\begin{equation}
{\theta}_j=-e^{-\tpi  \beta_j}, \quad{\phi}_j= -e^{\tpi  \alpha_j},\quad j=1,\ldots, g.
\label{eq:periodicities}
\end{equation}
It follows from (\ref{theta_char_periodicity}) that (\ref{Szegodefn}) is a function of $e^{\tpi  \alpha_i}$ and $e^{\tpi  \beta_i}$. The further factors of $-1$ in (\ref{eq:periodicities}) are included for later convenience.
The Szeg\"o kernel has multipliers along the $a_i$ and $b_j$ cycles in $x$ given by $-\phi_i$ and $-\theta_j$ respectively and is 
a meromorphic $(\half,\half)$-form on $\widetilde{\Sigma} \times \widetilde{\Sigma}$ satisfying:  
\begin{eqnarray}
  S\left[{\theta}\atop {\phi} \right] (x,y) 
&\sim &  \frac{1}{x-y}dx^{\half}dy^{\half}
\quad \mbox{for } x\sim y,
 \label{Sz_local}\\
 S\left[{\theta}\atop {\phi} \right] (x,y) 
& = &   -S\left[{\theta^{-1}}\atop {\phi^{-1}} \right] (y,x),
\label{Sz_skewsym}
\end{eqnarray}
where $\theta^{-1}=({\theta}_i^{-1})$ and $\phi^{-1}=(\phi_i^{-1})$. Note that the skew-symmetry property (\ref{Sz_skewsym}) implies $S\left[{\theta}\atop {\phi} \right] (x,y)$ has multipliers along the $a_i$ and $b_j$ cycles in $y$ given by $-\phi_i^{-1}$ and $-\theta_j^{-1}$ respectively. 

Finally, we describe the modular invariance of the Szeg\"o kernel under the symplectic group $Sp(2g,\mathbb{Z})$ where we find \cite{F1}
\begin{equation}
S\left[ {\tilde{\theta} \atop \tilde{\phi} }\right] (x, y\vert \tilde\Omega)=
S\left[ {{\theta} \atop {\phi} }\right] (x, y\vert \Omega),
\label{Szmod}
\end{equation}
with $\tilde\Omega$ of (\ref{eq:modOmega}) and  where $\tilde{\theta}_j=-e^{-\tpi  \tilde{\beta}_j}$, $\tilde{\phi}_j= -e^{\tpi  \tilde{\alpha}_j}$ for
\begin{equation}
\left( {-\tilde{\beta} \atop \tilde{\alpha} }\right)
=
\left(\begin{array}{cc}
A & B \\ 
C & D
\end{array}\right)
\left( {-\beta}\atop {\alpha} \right) 
+\half 
\left(\begin{array}{c}
{-\mbox{diag}(AB^T)} \\ 
{\mbox{diag}(CD^T)}
\end{array}\right),
\label{albetatilde}
\end{equation}
where $\mbox{diag}(M)$ denotes the diagonal elements of a matrix $M$.

\medskip
For a  Riemann surface of genus one described by an oriented torus $\mathbb{C}/{\Lambda}$ for lattice ${\Lambda}=\tpi (\mathbb{Z}\tau\oplus \mathbb{Z})$ for $\tau\in \mathbb{H}_{1}$, the genus one prime form is $E^{(1)}(x,y)=K(x-y,\tau )dx^{-\half}dy^{-\half}$
where 
\begin{equation}
K(z,\tau )=\frac{\vartheta _{1}(z,\tau )}{\partial_{z}\vartheta _{1}(0,\tau)},  
\label{Ktheta}
\end{equation}%
for $z\in \mathbb{C}$ and $\tau \in \mathbb{H}_{1}$ and where 
$\vartheta_{1}(z,\tau )
=\vartheta \left[{\half}\atop {\half} \right](z,\tau )$.

For $(\theta ,\phi )\neq (1,1)$ with $\theta =-e^{-\tpi \beta}$ and $\phi=-e^{\tpi \alpha}$ the genus   one  Szeg\"{o} kernel is
\begin{equation}
S^{(1)}\left[{\theta}\atop {\phi} \right] (x,y|\tau )=
P_{1}\left[{\theta}\atop {\phi} \right](x-y,\tau )dx^{\half}dy^{\half},\label{S1}
\end{equation}%
where
\begin{eqnarray}
P_{1}\left[{\theta}\atop {\phi} \right] (z,\tau )\nonumber
&=&\frac{\vartheta \left[{\alpha}\atop {\beta} \right](z,\tau )}
{\vartheta \left[{\alpha}\atop {\beta} \right](0,\tau )}
\frac{1}{K(z,\tau)},\nonumber\\
&=&-\sum\limits_{k\in\mathbb{Z}}\frac{q_z^{k+\lambda}}{1- \theta^{-1} q^{k+\lambda}},\label{P1}
\end{eqnarray}%
is a \lq twisted\rq\ Weierstrass function \cite{MTZ} for $q_z=e^z$ and with
$\phi =\exp (\tpi \lambda )$ for $0\leq \lambda <1$. 
The genus  one  modular group $SL(2,\mathbb{Z})$ acts in this case with (\ref{Szmod}) and (\ref{albetatilde}) following from 
\begin{equation}
P_{1}  \left(
\gamma\left[{\theta}\atop {\phi} \right]\right)(\gamma z|\gamma\tau )=(c\tau +d)
P_{1} \left[{\theta}\atop {\phi} \right](z|\tau ),
\label{eq:S1mod}
\end{equation}
with  
\begin{equation}
\gamma\tau=\frac{a\tau+b}{c\tau+d},\quad \gamma z=\frac{z}{c\tau+d},
\label{eq:modztau}
\end{equation} for  
$\gamma=\left(\begin{array}{cc}
a &  b \\ 
c  & d 
\end{array}%
\right)\in SL(2,\mathbb{Z})$ and 
\begin{equation}
\gamma\left[{\theta}\atop {\phi} \right]
=\left[{\theta^a\phi^b}\atop {\theta^c\phi^d} \right]. 
\label{modthetaphi}
\end{equation}

We also have a Laurant
expansion \cite{MTZ}
\begin{equation}
P_{1}\left[{\theta}\atop {\phi} \right](z,\tau )=\frac{1}{z}-\sum\limits_{n\geq 1}
E_{n}\left[{\theta}\atop {\phi} \right](\tau )z^{n-1},\label{P1zn}
\end{equation}
for twisted Eisenstein series defined  by  
\begin{eqnarray}
E_{n}\left[{\theta}\atop {\phi} \right](\tau ) &=&-\frac{B_{n}(\lambda )}{n!}+\frac{1}{(n-1)!}%
\sum\limits_{r\geq 0}\frac{(r+\lambda )^{n-1}\theta ^{-1}q^{r+\lambda }}{%
1-\theta ^{-1}q^{r+\lambda }} \notag\\
&&+\frac{(-1)^{n}}{(n-1)!}\sum\limits_{r\geq 1}\frac{(r-\lambda
)^{n-1}\theta q^{r-\lambda }}{1-\theta q^{r-\lambda }},\label{Endef}
\end{eqnarray}
for $n\geq 1$ and where $B_{n}(\lambda )$ is the Bernoulli polynomial defined by %%
\begin{equation*}
\frac{q_{z}^{\lambda }}{q_{z}-1}=\frac{1}{z}
+\sum\limits_{n\geq 1}\frac{B_{n}(\lambda )}{n!}z^{n-1}.
\end{equation*}
For $(\theta,\phi)=(1,1)$ and $n\ge 2$ the twisted Eisenstein series reduce to the standard elliptic Eisenstein series  with $E_n(\tau)=0$ for $n$ odd.

%%%%%%%%%%%%%%%%%%%%%%%%%%%%%%%%%%%%%%%%%%%%%%%%%%%%%%%%%%%%%%%%%%%%%%%
\section{The Szeg\"o Kernel on Two Sewn Riemann Surfaces}
%%%%%%%%%%%%%%%%%%%%%%%%%%%%%%%%%%%%%%%%%%%%%%%%%%%%%%%%%%%%%%%%%%%%%%%%
\label{Sect_Epsilon_g}
%%
%%%%%%%%%%%%%%%%%%%%%%%%%%%%%%%%%%%%%%%%%%%%%%%%%%%%%%%%%%%%%%%%%%%%%

%%
%%

%%%%%%%%%%%%%%%%%%%%%%%%%%%%%%%%%%%%%%%%%%%%%%%%%%%%%%%%%%%%%%%%%%%%%%%%%%%%%%
\subsection{The $\epsilon$-Formalism Sewing Scheme}
%%%
\label{Subsec_eps_sew}
%%%%%%%%%%%%%%%%%%%%%%%%%%%%%%%%%%%%%%%%%%%%%%%%%%%%%%%%%%%%%%%%%%%%%%%%%%%%%%
%%
We review the Yamada \cite{Y} formalism
 for \lq sewing\rq\ together two Riemann surfaces 
$\Sigma^{(g_a)}$ of genus $g_{a}$ for $a=1,2$ to form a surface of genus $g_{1}+g_{2}$. Following \cite{MT1}, we refer to this sewing scheme as the $\epsilon$-formalism. 

Choose a local coordinate $z_{a}$ on $\Sigma^{(g_a)}$ in the
neighborhood of a point $p_{a}$, and consider the closed disk $\left\vert
z_{a}\right\vert \leq r_{a}$ for $r_{a}>0$, sufficiently small. Let $\epsilon $ be a complex sewing
parameter with $|\epsilon |\leq r_{1}r_{2}$ and excise the disk 
\begin{equation*}
\{z_{a}:\, \left\vert z_{a}\right\vert \leq |\epsilon |/r_{\bar{a}}\}\subset 
\Sigma^{(g_a)}, 
\end{equation*}
to form a punctured surface 
\begin{equation*}
\widehat{\Sigma}^{(g_a)}=\Sigma^{(g_a)}\backslash \{z_{a}:\, \left\vert
z_{a}\right\vert \leq |\epsilon |/r_{\bar{a}}\}.
\end{equation*}
Here and below, we use the convention 
\begin{equation*}
\overline{1}=2,\quad \overline{2}=1.  
\label{bardef}
\end{equation*}
Define the annulus $\mathcal{A}_{a}=\{z_{a}:\, |\epsilon |/r_{\bar{a}}\leq \left\vert
z_{a}\right\vert \leq r_{a}\}\subset \widehat{\Sigma}^{(g_a)}
$
and identify $\mathcal{A}_{1}$ and $\mathcal{A}_{2}$ as a single region $%
\mathcal{A}=\mathcal{A}_{1}\simeq \mathcal{A}_{2}$ via the sewing relation 
\begin{equation}
z_{1}z_{2}=\epsilon .  \label{pinch}
\end{equation}
\medskip

\begin{center}
\begin{picture}(300,100)

%left surface
\put(50,52){\qbezier(-30,18)(-10,10)(10,18)}% left left upper
\put(50,52){\qbezier(10,18)(50,35)(90,18)}% left upper
\put(50,48){\qbezier(-30,-18)(-10,-10)(10,-18)}%left left lower
\put(50,48){\qbezier(10,-18)(50,-35)(90,-18)}%left lower

\put(45,50){\qbezier(25,0)(45,17)(60,0)}%upper
\put(45,50){\qbezier(20,2)(45,-17)(65,2)}%lower

%right surface
\put(175,52){\qbezier(90,18)(110,10)(130,18)}% right right upper
\put(175,52){\qbezier(10,18)(50,35)(90,18)}%right upper
\put(175,48){\qbezier(90,-18)(110,-10)(130,-18)}%right right lower
\put(175,48){\qbezier(10,-18)(50,-35)(90,-18)}% right lower

\put(200,50){\qbezier(25,0)(45,17)(60,0)}%upper
\put(200,50){\qbezier(20,2)(45,-17)(65,2)}%lower
% left annulus centered at (140,50)

\put(140,50){\circle{16}}
\put(140,50){\circle{40}}

% z_{1} =0label
\put(140,50){\vector(-1,-2){0}}%arrow
\put(50,50){\qbezier(90,0)(100,15)(90,30)}%
\put(140,90){\makebox(0,0){$z_{1} =0$}}

% line and r1 label
\put(140,50){\line(-1,1){14.1}}
\put(127,55){\makebox(0,0){$r_{1} $}}

% line and eps/r2 label
\put(140,50){\line(1,0){8}}
\put(145,50){\vector(1,4){0}}%arrow
\put(55,20){\qbezier(90,4)(85,17)(90,30)}%
\put(150,15){\makebox(0,0){$|\epsilon|/r_2$}}

%contour C1
%\put(140,50){\circle{28}}
%\put(154,50){\vector(0,1){0}}%arrow
%\put(154,55){\makebox(0,0){$\mathcal{C}_{1} $}}

%Sg1 label
\put(50,50){\makebox(0,0){$\widehat{\Sigma}^{(g_{1} )}$}}

%right annulus centred at (185,50)

\put(185,50){\circle{16}}
\put(185,50){\circle{40}}

% z_2=0label
\put(185,50){\vector(1,-2){0}}%arrow
\put(95,50){\qbezier(90,0)(80,15)(90,30)}%
\put(185,90){\makebox(0,0){$z_2=0$}}

% line and r2 label
\put(185,50){\line(-1,-1){14.1}}
\put(171,45){\makebox(0,0){$r_2$}}

% line and eps/r1 label
\put(185,50){\line(1,0){8}}
\put(190,50){\vector(1,4){0}}%arrow
\put(90,20){\qbezier(100,4)(95,17)(100,30)}%
\put(190,15){\makebox(0,0){$|\epsilon|/r_{1} $}}

%Sg2 label
\put(280,50){\makebox(0,0){$\widehat{\Sigma}^{(g_2)}$}}

\end{picture}

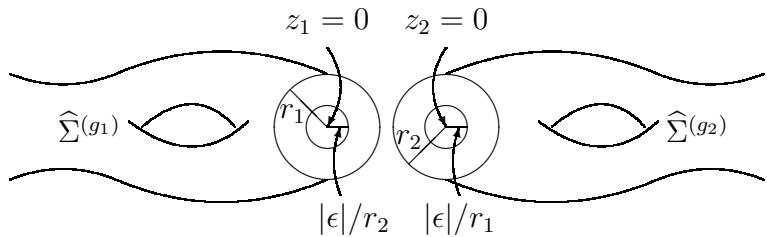
{\small Fig. 1: Sewing Two Riemann Surfaces}
\end{center}
In this way we obtain a compact Riemann surface $\Sigma^{(g_{1} +g_2)}=\{ \widehat {\Sigma}^{g_{1} }
\backslash \mathcal{A}_{1}\}
\cup \{\widehat{\Sigma}^{(g_2)}\backslash 
\mathcal{A}_{2}\}\cup \mathcal{A}$ of genus $g_{1}+g_{2}$. By construction, $\Sigma^{(g_{1} +g_2)}$ degenerates into $\Sigma^{(g_{1})}$ and $\Sigma^{(g_2)}$ in the limit $\epsilon\rightarrow 0$.

\medskip

The form $\omega^{(g_{1} +g_2)}$ on $\Sigma^{(g_{1} +g_2)}$ can be found in terms of data coming from $\omega^{(g_a)}$ on $\widehat{\Sigma}^{(g_a)}$ \cite{Y}. $\Sigma^{(g_{1} +g_2)}$ inherits a homology cycle basis labeled $\{a_{s_{1}},b_{s_{1}}:s_{1}=1,\ldots, g_{1}\}$ and 
$\{a_{s_{2}},b_{s_{2}}:s_{2}=g_{1}+1,\ldots, g_{1}+g_{2}\}$ from $\Sigma^{(g_{1} )}$ and $\Sigma^{(g_2)}$ respectively. This allows us to compute the normalized 1-forms $\nu_{i}^{(g_{1} +g_2)}$ and the period matrix $\Omega_{ij}^{(g_{1} +g_2)}$. In particular, we find \cite{Y},\cite{MT1}
%%%%%%%%%%%%%%%%%%%%%%%%%%%%%%%%%%%%%%%%%%%%%%%%%
\begin{theorem}
\label{theoremomg1g2om} 
$\omega^{(g_{1}+g_{2})}$, $\nu_{i}^{(g_{1} +g_2)}$ and $\Omega_{ij}^{(g_{1} +g_2)}$ are holomorphic in $\epsilon $ for 
$|\epsilon |<r_{1}r_{2}$ with
\begin{eqnarray*}
&&\omega^{(g_{1}+g_{2})}(x,y)=
\delta_{ab}\omega^{(g_{a})}(x,y)+O(\epsilon),\\
&&\nu^{(g_{1}+g_{2})}_{s_{b}}(x)=
\delta_{ab}\nu^{(g_{a})}_{s_{a}}(x)+O(\epsilon),\\
&&\Omega^{(g_{1}+g_{2})}_{s_{a}t_{b}}=
\delta_{ab}\Omega^{(g_{a})}_{s_{a}t_{a}}+O(\epsilon),
\end{eqnarray*}  
for $x\in \widehat{\Sigma}^{(g_a)}$, $y\in \widehat{\Sigma}^{(g_b)}$ and $a,b=1,2$ and
where $s_{a},t_{b}$ label the inherited homology basis.   
\end{theorem} 
The explicit form of $\omega^{(g_{1} +g_2)}$, $\nu_{i}^{(g_{1} +g_2)}$ and $\Omega_{ij}^{(g_{1} +g_2)}$ is described in \cite{Y},\cite{MT1}. 
%%
%%%%%%%%%%%%%%%%%%%%%%%%%%%%%%%%%%%%%%%%%%%%%%%%%%%%%%%%%%%%%%%%%%%%%%%%%%%
%%
%%%%%%%%%%%%%%%%%%%%%%%%%%%%%%%%%%%%%%%%%%%%%%%%%%%%%%%%%%%%%%%%%%%%%%%%%%%%%%%%%%%%
\subsection{The Szeg\"o Kernel in the $\epsilon$-Formalism}  
%%%%%%%%%%%%%%%%%%%%%%%%%%%%%%%%%%%%%%%%%%%%%%%%%%%%%%%%%%%%%%%%%%%%%%%%%%%%%%%%%%%%%%
%%
\label{Subsec_eps_Szego}
We now determine the Szeg\"o kernel on the Riemann surface $\Sigma^{(g_{1} +g_2)}$ in terms of data coming from Szeg\"o kernel 
\begin{eqnarray}
S^{(g_a)}(x,y)=S^{(g_a)}\left[ {\theta}^{(g_a)}  \atop {\phi}^{(g_a)} \right](x,y),
\label{eq:Szego_a}
\end{eqnarray}
on the surface $\Sigma^{(g_a)}$ for $a=1,2$.  
We adopt the abbreviated notation of the left hand side of (\ref{eq:Szego_a}) when there is no ambiguity. Similarly, the Szeg\"o kernel on $\Sigma^{(g_{1} +g_2)}$ is denoted by
\begin{eqnarray}
S^{(g_{1} +g_2)}(x,y)=
S^{(g_{1}  + g_2)}\left[ {\theta}^{(g_{1} +g_2)}  \atop {\phi}^{(g_{1} +g_2)} \right] (x,y),
\label{eq:Szego_g12}
\end{eqnarray}
with periodicities $({\theta}^{(g_{1} +g_2)}_{s_a}, {\phi}^{(g_{1} +g_2)}_{s_a})= ({\theta}^{(g_a)}_{s_a}, {\phi}^{(g_a)}_{s_a})$ on the inherited homology basis.\footnote{Note that we exclude those Riemann theta characteristics for which (\ref{eq:Szego_g12}) exists but where either of the lower genus theta functions vanishes i.e. we assume that (\ref{eq:Szego_a}) exists for $a=1,2$.} 
\medskip

We next describe $S^{(g_{1} +g_2)}(x,y)$ in terms of $S^{(g_a)}(x,y)$. We first show that 
%%%%%%%%%%%%%%%%%%%%%%%%%%%%%%%%%%%%%%%%%%%%%%%%%
\begin{theorem}
\label{theoremomg1g2holo} 
$S^{(g_{1}+g_{2})}$ is holomorphic in $\epsilon^{\half} $ for 
$|\epsilon |<r_{1}r_{2}$ with
\begin{equation*}
S^{(g_{1}+g_{2})}(x,y)=\left\{
\begin{array}{l}
S^{(g_{a})}(x,y)+O(\epsilon),\quad \mbox{ for } x,y\in \widehat{\Sigma}^{(g_a)}, 
\\
O(\epsilon^{\half}), \quad \mbox{ for }
 x\in \widehat{\Sigma}^{(g_a)},
\; y\in \widehat{\Sigma}^{(g_{\bar{a}})}. 
\end{array}\right.
\end{equation*}   
\end{theorem}
%%%%%%%%%%%%%%%%%

\noindent 
{\bf Proof.}
Applying Theorem \ref{theoremomg1g2om} to (\ref{theta}) we have
  \begin{eqnarray}
\vartheta \left[ {{\alpha^{(g_{1} +g_2)}} \atop {\beta^{(g_{1} +g_2)}} }\right]  \left( {z}^{(g_{1} +g_2)}\vert  \Omega^{(g_{1} +g_2)}  \right)&=&
\vartheta \left[ {{\alpha^{(g_{1} )}} \atop {\beta^{(g_{1} )}} }\right]  \left( {z}^{(g_{1} )}\vert  \Omega^{(g_{1} )}  \right)
\vartheta \left[ {{\alpha^{(g_2)}} \atop {\beta^{(g_2)}} }\right]  \left( {z}^{(g_2)}\vert  \Omega^{(g_2)}  \right)\notag\\
&&+O(\epsilon),
\label{theta_eps}  
\end{eqnarray}
with $({\alpha}^{(g_{1} +g_2)})=(\alpha^{(g_{1} )}_{1}, \ldots, \alpha^{(g_{1} )}_{g_{1} },\alpha^{(g_2)}_{1}, \ldots, \alpha^{(g_2)}_{g_2})$ etc.
We firstly show that the genus $g_{1} +g_2$ prime form obeys
\begin{equation}
E^{(g_{1} +g_2)}(x,y)
=\left\{ 
\begin{array}{l}
E^{(g_a)}(x,y)+O(\epsilon),  \quad \mbox{ for } x,y \in 
\widehat{\Sigma}^{(g_a)}, 
\\
O(\epsilon^{-{\half} }), \quad \mbox{ for }
 x\in \widehat{\Sigma}^{(g_a)},
\; y\in \widehat{\Sigma}^{(g_{\bar{a}})}. 
\end{array} 
\right. \label{prime_eps}
\end{equation}
For the genus $g_{1} +g_2$ odd characteristic of (\ref{eq:thetadelta}) we find from 
(\ref{theta_eps}) that either 
$\vartheta \left[ { \gamma}^{(g_{1} )} \atop {\delta}^{(g_{1} )}  \right](0)\neq 0$ or
$\vartheta \left[ { \gamma}^{(g_2)} \atop {\delta}^{(g_2)}  \right](0)\neq 0$ on the lower genus surfaces. Hence it follows that $\zeta^{(g_{1} +g_2)}(x)\zeta^{(g_{1} +g_2)}(y)=O(\epsilon)$ for $x\in \widehat{\Sigma}^{(g_a)},
\; y\in \widehat{\Sigma}^{(g_{\bar{a}})}$ for the 1-form (\ref{eq:zeta}). 
We also note that 
\begin{equation}
\int_y^x\nu_{s_b}^{(g_{1} +g_2)}
=\left\{ 
\begin{array}{l}
\int_y^x\nu_{s_a}^{(g_a)}+O(\epsilon),  \quad \mbox{ for }  x,y \in 
\widehat{\Sigma}^{(g_a)}, 
\\
\delta_{ab}\int_{p_a}^x\nu_{s_a}^{(g_a)}+\delta_{{\bar{a}}b}\int^{p_{\bar{a}}}_y\nu_{s_{\bar{a}}}^{(g_{\bar{a}})}+O(\epsilon), \quad \mbox{ for }
 x\in \widehat{\Sigma}^{(g_a)},
\; y\in \widehat{\Sigma}^{(g_{\bar{a}})}, 
\end{array} 
\right. \label{nu_eps}
\end{equation}   
from which it  follows  that $E^{(g_{1} +g_2)}(x,y)=O(\epsilon^{-\half})$ for 
$x\in \widehat{\Sigma}^{(g_a)}$ and 
$y\in \widehat{\Sigma}^{(g_{\bar{a}})}$. 

We next determine $E^{(g_{1} +g_2)}(x,y)$  for  $x,y \in \widehat{\Sigma}^{(g_a)}$. The differential $\omega^{(g_{1}+g_{2})}$ for 
$x,y \in \widehat{\Sigma}^{(g_a)}$ obeys 
\begin{equation}
\omega^{(g_{1}+g_{2})}(x,y)-\omega^{(g_a)}(x,y)=
a_a(x)X_{{\bar{a}}{\bar{a}}}a^T_a(y)=O(\epsilon),\label{om_eps}
\end{equation}
where $a_a(x)X_{{\bar{a}}{\bar{a}}}a^T_a(y)=\sum_{k,l\ge 1}a_a(x,k)X_{{\bar{a}}{\bar{a}}}(k,l)a_a(y,l)$ with  $a_a(x,k)$ a certain 1-form on $\widehat{\Sigma}^{(g_a)}$
and $X_{{\bar{a}}{\bar{a}}}(k,l)$ an infinite matrix determined from genus $g_{1} $ and $g_2$ data (see \cite{MT1} for details). It follows from (\ref{eq:Prime_omega}) that
\begin{equation*}
E^{(g_{1}+g_{2})}(x,y)=E^{(g_a)}(x,y) 
e^{-\half b_aX_{{\bar{a}}{\bar{a}}}b_a^T}=E^{(g_a)}(x,y)+O(\epsilon),
\end{equation*} 
where $b_a(k)=\int_y^x a_a(\cdot,k)$.
Thus (\ref{prime_eps}) holds.  We then apply Theorem \ref{theoremomg1g2om},  (\ref{theta_eps}),  (\ref{prime_eps}) and (\ref{nu_eps}) to  (\ref{Szegodefn}) to prove the result.
\hfill $\square$

%%%%%%%%%%%%%%%%%%%%%%%%%%%%%%%%%%%%%%%%%%
\medskip
We next remark that for $x,z_a\in \widehat{\Sigma}^{(g_a)}$ then $S^{(g_{a})}(x,z_a)S^{(g_{1}+g_2)}(z_a, y)$ is a meromorphic 1-form (cf. \cite{HS}) in $z_a$ periodic on the $\Sigma^{(g_a)}$ cycles (cf. (\ref{Sz_skewsym})) with simple poles described by (\ref{Sz_local}) where
\begin{eqnarray}
S^{(g_{a})}(x,z_a)S^{(g_{1}+g_2)}(z_a,y)&\sim& \frac{dz_a}{x-z_a}S^{(g_{1}+g_2)}(x,y)
\mbox{ for }z_a \sim x,\notag\\
S^{(g_{a})}(x,z_a)S^{(g_{1}+g_2)}(z_a,y)&\sim& \frac{dz_a}{z_a-y}S^{(g_{a})}(x,y) 
\mbox{ for $z_a \sim y$ if $ y\in \widehat{\Sigma}^{(g_a)}$}. \notag \\
\label{eq:SSpoles}
\end{eqnarray}   
A similar behavior holds for $S^{(g_{1}+g_2)}(x,z_b)S^{(g_{b})}(z_b,y)$ as  a meromorphic 1-form in $z_b$. This allows us to determine the following integral equations
%%%%%%%%%%%%%%%%%%%%%%%%%%%
\begin{proposition}
%%%%%%%%%%%%%%%%%%%%%%%%%%%%%%
\label{Prop_inteqn1} 
The Szeg\"o kernel on $\Sigma^{(g_{1} + g_2)}$ is given by  
\begin{eqnarray}
S^{(g_{1}+g_{2})}(x,y)&=&\delta _{ab}S^{(g_{a})}(x,y) 
- \frac{1}{\tpi }\oint\limits_{\mathcal{C}_{a}(z_{a})} 
S^{(g_{a})}(x,z_a) 
S^{(g_{1}+g_2)}(z_a, y),\label{Prop_epsinteq1} 
\\
&=& \delta _{ab}S^{(g_{a})} (x,y)   
+\frac{1}{\tpi }\oint\limits_{\mathcal{C}_{b}(z_{b})} 
S^{(g_{1}+g_{2})}(x,z_b) 
S^{(g_{b})}( z_b, y),\label{Prop_epsinteq2}  
\end{eqnarray}
for $x\in \widehat{\Sigma}^{(g_a)}$, $y\in \widehat{\Sigma}^{(g_b)}$ for $a,b=1,2$ and where $\mathcal{C}_{a}(z_{a})\subset \mathcal{A}_{a}$ denotes a closed
anti-clockwise oriented contour parameterized by $z_{a}$ surrounding the
puncture at $z_{a}=0$ on $\widehat{\Sigma}^{(g_a)}$.  
\end{proposition}
\noindent 
{\bf Proof.}  
Let $\sigma_a$ be a contour on $\widehat{\Sigma}^{(g_a)}$ surrounding $\mathcal{A}_{a}$ and the given points $x$ (and $y$, if $a=b$) on $\widehat{\Sigma}^{(g_a)}$ (see Fig.~2).\footnote{$\sigma_a$ may be construed as being the boundary of the simple-connected covering space for $\Sigma^{(g_a)}$ as illustrated in Fig.~2 for a genus two surface. }    
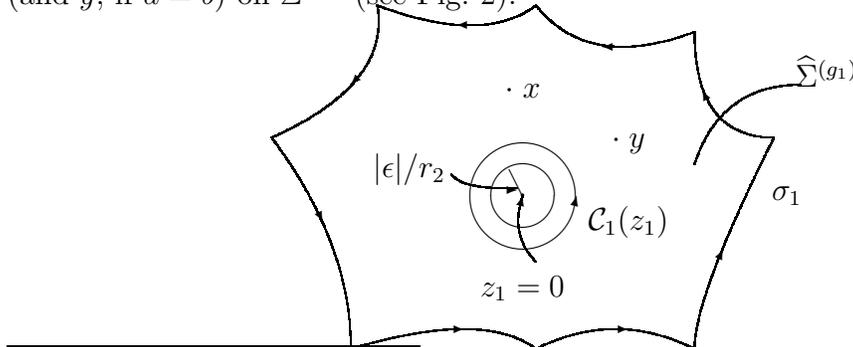
\medskip

%%%%%%%%%%%%%%%%%%%%%%%%%%%%%%%%%%%%%%%%%%%%%%%%%%%%%%%%%%%%%%%%%%%%%% 
\begin{center}

\medskip
\begin{picture}(300,100)
%contour figure

\put(175,52){\qbezier(-120,20)(-75, 40)(-80,70)}%

\put(175,52){\qbezier(-80,70)(-40, 55)(-20,70)}%

\put(175,52){\qbezier(-20,70)(0, 45)(40,60)}%

\put(175,52){\qbezier(40,60)(40, 20)(70,20)}%

%%%%%%%%%%%%%%%%%%%%%%%%%%%%%%%%%%%%%%%%%%%%%%%%%%%%%%%%
\put(175,52){\qbezier(-120,20)(-90, -20)(-90,-60)}%

\put(175,52){\qbezier(-90,-60)(-40, -45)(-20,-60)}%

\put(175,52){\qbezier(-20,-60)(10, -45)(40,-60)}%

\put(175,52){\qbezier(40,-60)(40, -40)(70, 20)}%
%%%%%%%%%%%%%%%%%%%%%%%%%%%%%%%%%%%%%%%%%%%%%

\put(175,52){\qbezier(40,10)(50, 40)(80, 40)}%

\put(175,52){\makebox(180,90){${\widehat \Sigma}^{(g_{1} )}$}}%

%%\put(235,84){\vector(-1,-1){0}}%arrow

%%%%%%%%%%%%%%%%%%%%%%%%%%%%%%%%%%%%%%%%%%%%%
% left annulus 

%%0
\put(150,50){\circle{40}}
\put(150,50){\circle{22}}

\put(170.5,50){\vector(1,4){0}}%arrow

%%0
\put(150,50){\makebox(0,0){$\cdot$}}

\put(150,15){\makebox(0,0){${z_{1} =0}$}}

%%0
\put(150,50){\vector(1,4){0}}%arrow

%%0
{\qbezier(155,25)(145,35)(150,50)}%

\put(150,90){\makebox(0,0){$\cdot \; x$}}

%%%%%%%%%%%%%%%%%%%%%%%%%%%%%%%%%%%%%%%%%%%
\put(85,92){\vector(-1, -1){0}}%arrow

\put(125,114.5){\vector(-1, 0){0}}%arrow

\put(180,106.5){\vector(-1, 0){0}}%arrow

\put(218,90){\vector(-1, 2){0}}%arrow
%%

%%%%%%%%%%%%%%%%%%%%%%%%%%%%%%%%%%%%%%%%%%%
\put(128,-0.5){\vector(1, 0){0}}%arrow

\put(190,-0.5){\vector(1, 0){0}}%arrow

\put(74.5,40){\vector(1, -2){0}}%arrow

% line and r2 label

%%0
\put(190,40){\makebox(0,0){$\mathcal{C}_{1}(z_{1} )$}}

\put(190,70){\makebox(0,0){$\cdot \; y$}}

\put(250,50){\makebox(0,0){$\sigma_{1} $}}

\put(226,30){\vector(1,2){0}}%arrow

%%%%%%%%%%%%%%%%%%%%%%%%%%%%%%%%%%%%%%%%%%%%%%%%%%%
% line and rho/r2 label
%%0
\put(150,50){\line(-1,2){5}}

%%0
\put(148,52){\vector(1,0){0}}%arrow

%%0
\put(43,33){\qbezier(105, 20)(85,19)(80,25)}%

%%0
\put(107,60){\makebox(0,0){$|\epsilon|/r_{2}$}}

\end{picture}
\vskip 0.5 truecm
\small{Fig. 2: Example with $x,y\in \widehat\Sigma^{(g_{1} )}$}
\end{center}
From Cauchy's theorem $\oint_{\sigma_a}
S^{(g_{a})}(x,z_a)
S^{(g_{1}+g_2)}(z_a, y) = 0$ and hence \eqref{eq:SSpoles} gives 
\begin{eqnarray*}
0=-S^{(g_{1}+g_2)}(x,y)+\delta_{ab}S^{(g_a)}(x,y)+\frac{1} {\tpi }\oint_{\mathcal{C}_a(z_a)} 
S^{(g_{a})}(x,z_a)
S^{(g_{1}+g_2)}(z_a, y),
\end{eqnarray*}
giving (\ref{Prop_epsinteq1}). Considering $S^{(g_{1}+g_2)}(x,z_b)S^{(g_{b})}(z_b,y)$ leads to 
(\ref{Prop_epsinteq2}).   
\hfill  $\square$ 

%%%%%%%%%%%%%%%%%%%%%%%%%%%%%%%%%%%%%%%%%%%%%%%%%%%%%%
\medskip
Similarly to \cite{MT1} we define weighted moments for $S^{(g_{1}+g_{2})}$ by
\begin{eqnarray}
 X_{ab}(k,l,\epsilon)=&& X_{ab}\left[ {\theta}^{(g_{1} +g_2)}  \atop {\phi}^{(g_{1} +g_2)} \right](k,l,\epsilon)\notag \\
=&&  \frac{\epsilon^{\half(k+l-1)} }  { (\tpi )^{2} }  
\oint\limits_{ \mathcal{C}_{a}(x)}\oint\limits_{\mathcal{C}_{b}(y)}x^{-k}y^{-l}
S^{(g_{1}+g_{2})}(x,y) dx^{\half} dy^{\half},  
\label{Xijdef}
\end{eqnarray} 
%%%%%%%%%%%%%%%%%%%%%%%%%%%%%%%%%%%%%%%%%%%%%%%%%%%%%%%%
for $k,l \ge 1$. From (\ref{Sz_skewsym}) it follows that 
\begin{equation}
X_{ab}\left[ {\theta}^{(g_{1} +g_2)}  \atop {\phi}^{(g_{1} +g_2)} \right](k,l,\epsilon)=
-X_{ba}\left[ ({\theta}^{(g_{1} +g_2)})^{-1}  \atop ({\phi}^{(g_{1} +g_2)^{-1}}) \right](l,k,\epsilon).
\label{Xabsym}
\end{equation}
We denote by $X_{ab}=(X_{ab}(k,l, \epsilon))$ the infinite matrix indexed by $k,l\ge 1$. 

We also define various moments for $S^{(g_{a})}(x,y)$. These provide the data used to construct $S^{(g_{1}+g_2)}(x,y)$.  Define holomorphic $\half$-forms on $\widehat{\Sigma}^{(g_a)}$  by
\begin{eqnarray}
 h_{a}(k,x,\epsilon )=&& h_{a}\left[ {\theta}^{(g_a)}  \atop {\phi}^{(g_a)} \right](k,x,\epsilon)
 =\frac{ \epsilon^{\frac{k}{2} -\frac{1}{4}} }{\tpi } 
\oint\limits_{
\mathcal{C}_{a}(z_{a})}S^{(g_{a})}(x,z_{a}) z_{a}^{-k} dz_a^{\half}, 
 \label{hdef}\\
 \bar{h}_{a}(k,y,\epsilon )=&& \bar{h}_{a}\left[ {\theta}^{(g_a)}  \atop {\phi}^{(g_a)} \right](k,y,\epsilon)=
 \frac{  \epsilon^{\frac{k}{2} -\frac{1}{4}}}{\tpi } 
\oint\limits_{
 \mathcal{C}_{a}(z_{a})}S^{(g_{a})}(z_{a}, y) z_{a}^{-k}  dz_a^{\half}, 
 \label{barhdef}
\end{eqnarray}
%%%%%%%%%%%%%%%%%%%%%%%%%%%%%%%%%%%%%%%
and introduce infinite row vectors $h_{a}(x)=(h_{a}(k,x))$,  
$\bar{h}_{a}(x)=(\bar{h}_{a}(k,x))$ indexed by $k\ge 1$. From (\ref{Sz_skewsym}) it follows that 
\begin{equation}
 \bar{h}_a\left[ {\theta}^{(g_a)}  \atop {\phi}^{(g_a)}\right](k,x,\epsilon) = 
-h_a\left[ ({\theta}^{(g_a)})^{-1}  \atop ({\phi}^{(g_a)})^{-1}\right](k,x,\epsilon).
\label{hhbar}
\end{equation}
Finally, we define the moment matrix  
\begin{eqnarray}
F_{a}( k,l,\epsilon) =&& F_{a}\left[ {\theta}^{(g_a)}  \atop {\phi}^{(g_a)} \right](k,l,\epsilon)\notag\\
=&&    \frac{\epsilon^{\half(k+l-1)}}{(\tpi )^{2}}  
\oint\limits_{\mathcal{C}_{a}(x)}\oint\limits_{\mathcal{C}_{a}(y)}x^{-k}y^{-l}
S^{(g_{a})}(x,y) dx^{\half} dy^{\half}
  \notag \\=&&      \frac{\epsilon^{\frac{k}{2} -\frac{1}{4}}}{\tpi  } 
\oint\limits_{\mathcal{C}_{a}(x)}x^{-k}h_{a}(l,x) dx^{\half}
%\notag\\
%=&&     
=\frac{\epsilon^{\frac{l}{2} -\frac{1}{4}} }{\tpi  } 
\oint\limits_{\mathcal{C}_{a}(y)}y^{-l}\bar{h}_{a}(k,y) dy^{\half}.\  
\label{Fdef}
\end{eqnarray}
$F_{a}( k,l,\epsilon)$ obeys a skew-symmetry property from (\ref{Sz_skewsym}) similar to (\ref{Xabsym}).
%%
%%%%%%%%%%%%%%%%%%%%%%%%%%%%%%%%%%
We may invert (\ref{hdef})--(\ref{Fdef}) using (\ref{Sz_local}) to find for  $x,y\in\widehat{\Sigma}^{(g_a)}$ that
\begin{eqnarray}
S^{(g_{a})}(x,y) &=&\left[\frac{1}{x-y} 
+\sum_{k,l\geq 1}\epsilon ^{-\half(k+l-1)}
F_{a}(k,l,\epsilon) x^{k-1}y^{l-1} \right]dx^{\half}dy^{\half}\label{FtoS}\\
&=& \sum_{k\geq 1} 
  \epsilon^{-\frac{k}{2} +\frac{1}{4}} h_{a}(k,x)y^{k-1}dy^{\half}\label{htoS}\\
&=& \sum_{l\geq 1} 
  \epsilon^{-\frac{l}{2} +\frac{1}{4}}x^{l-1}  \bar{h}_{a}(l,y) dx^{\half}\label{hbartoS}.  
\end{eqnarray}
\medskip

We are now in a position to express $S^{(g_{1}+g_2)}(x,y)$ in terms of the lower genus data. From the sewing relation (\ref{pinch}) we have $dz_a=-\epsilon \frac{dz_{\bar a}}{z_{\bar a}^2}$ so that
\begin{equation}
dz_a^{\half}=(-1)^{\bar a}\xi\epsilon^{\half} \frac{dz_{\bar a}^{\half}}{z_{\bar a}},
\label{dz1dz2}
\end{equation}
where $\xi\in\{\pm \sqrt{-1}\}$ determines the square root branch chosen. We then find
\begin{proposition}
\label{Prop_hXh} 
$S^{(g_{1}+g_{2})}(x,y)$ is given by 
\begin{equation}
\label{Sg1g2}
S^{(g_{1}+g_{2})}(x,y)
=\left\{ 
\begin{array}{l}
S^{(g_{a})}(x,y) + h_{a}(x) X_{\bar{a}\bar{a}}  \bar{h}_{a}^{T}(y),  \quad \mbox{ for } x,y \in 
\widehat{\Sigma}^{(g_a)}, 
\\
  h_{a}(x)\left(\xi (-1)^{\bar{a}}  I- X_{\bar{a}a}\right) \bar{h}_{\bar{a}}^{T}(y), \quad \mbox{ for }
 x\in \widehat{\Sigma}^{(g_a)},
\; y\in \widehat{\Sigma}^{(g_{\bar{a}})}, 
\end{array} 
\right. 
\end{equation}
where $I$ denotes the infinite identity matrix and $T$ the transpose.
\end{proposition}
%%
%%
%%%%%%%%%%%%%%%%%%%%%%%%%%%%%%%
\noindent \textbf{Proof.} 
Consider $x,y\in \widehat{\Sigma}^{(g_{1} )}$. Noting that $\mathcal{C}_{1}(z_{1})$ may be deformed to $-\mathcal{C}_{2}(z_{2})$  on $\cal A$ via (\ref{pinch}) we find
%%%%%%%%%%%%%%%%%%%%%%%%%%%%%%%%%%%%%%%%%%%%
\begin{eqnarray*}
&& S ^{(g_{1}+g_{2})}(x,y)-S^{(g_{1})}(x,y) 
= -  \frac{1}{\tpi }\oint\limits_{\mathcal{C}_{1}(z_{1})} S^{(g_{1})}(x,z_{1} ) 
S^{(g_{1} + g_{2})}(z_{1} , y)
\\
&=&  
-  \sum\limits_{ k \ge 1} 
h_{1} (k,x)
\frac{ \epsilon^{-\frac{k}{2} + \frac{1}{4}}}{\tpi } \oint\limits_{\mathcal{C}_{1}(z_{1})} S^{(g_{1}+g_{2})}(z_{1} , y) 
 z_{1} ^{k-1}  dz_{1} ^{\half}
\\
&=&
\xi 
\sum\limits_{ k \ge 1}
h_{1} (k,x)  
\frac{
\epsilon^{\frac{k}{2} - \frac{1}{4}}} {\tpi  } 
\oint\limits_{\mathcal{C}_{2}(z_{2})} 
S^{(g_{1}+g_{2})}(z_2, y)  z_2^{-k}  dz_2^{\half} 
\\
&=& \xi
\sum\limits_{ k \ge 1} 
 h_{1} (k,x)
\frac{ \epsilon^{\frac{k}{2} - \frac{1}{4}}}{(\tpi )^2}
 \oint\limits_{\mathcal{C}_{2}(z_{2})} 
\oint\limits_{\mathcal{C}_{1}(u_{1})}
S^{(g_{1}+g_{2})}( z_2, u_{1} )   S^{(g_{1})}( u_{1} , y)  
 z_2^{-k} dz_2^{\half}  
\\
&=& -\xi^2 \sum\limits_{ k,l \ge 1} 
 h_{1} (k,x)  \bar{h}_{1} (l,y) 
 \frac{\epsilon^{ \half(k+l-1) }}{(\tpi )^2} \oint\limits_{\mathcal{C}_{2}(z_{2})}
\oint\limits_{\mathcal{C}_{2}(u_{2})} 
  S^{(g_{1}+g_{2})}(z_2, u_2)  
 {z_2^{-k}} {u_2^{-l} }  {dz_2^{\half}}{du_2^{\half}}   
\\
&=& h_{1}(x)X_{22} \bar{h}_{1}^{T}(y), 
\end{eqnarray*}
%%%%%%%%%%%%%%%%%%%%%%%%%%%%%%%%%%%%%%%%%%%%%%%%
using (\ref{Prop_epsinteq1}), (\ref{htoS}), (\ref{dz1dz2}), (\ref{hbartoS}), (\ref{Prop_epsinteq2}) and (\ref{dz1dz2}) again, respectively. Thus we recover the first line of (\ref{Sg1g2}) for $a=b=1$. A similar analysis holds for $a=b=2$.
%%
%%%%%%%%%%%%%%%%%%%%%%%%%%%%%%%%%%%%%%%%%%%%%%%%

For $x\in \widehat{\Sigma}^{(g_{1} )}$, $y\in \widehat{\Sigma}^{(g_2)}$ we find that 
\begin{eqnarray*}
&& S^{(g_{1}+g_{2})}(x,y) = 
- \frac{1}{\tpi } \oint\limits_{\mathcal{C}_{1}(z_{1})} S^{(g_{1})}(x, z_{1} )
S^{(g_{1}+g_{2})}(z_{1} ,y) 
\\
&=&
%% 
%%%%  -
\xi\sum\limits_{k \ge 1}  
\frac{\epsilon^{\frac{k}{2} - \frac{1}{4}}}{\tpi  }\oint\limits_{\mathcal{C}_{2}(z_{2})}
  z_2^{-k} dz_2^{\half} h_{1} (k,x) 
\\
&& \qquad \qquad  
\cdot \left( S^{(g_{2})}(z_2,y) 
+ \frac{1}{\tpi }\oint\limits_{\mathcal{C}_{2}(u_2)} 
S^{(g_{1}  + g_{2})}(z_2, u_2)  S^{(g_{1})}(u_2,y)  
 \right) 
\\
&=&
 \xi 
%%%
\sum\limits_{k \ge 1} h_{1} (k,x)\bar{h}_2(k,y) 
\\
&& +\xi^2  
  \sum\limits_{k, l \ge 1} h_{1} (k,x) \bar{h}_2(l,y) 
\frac{\epsilon^{\half(k+l -1)}}{(\tpi )^2} 
\oint\limits_{\mathcal{C}_{2}(z_{2})}
\oint\limits_{\mathcal{C}_{1}(u_{1})}
S^{(g_{1}  + g_{2})}(z_2, u_{1} )  
z_2^{-k} dz_2^{\half}  u_{1} ^{-l} du_{1} ^{\half}
\\
&=&   h_{1} (x)\left( \xi I - X_{21}\right) \bar{h}_2^T(y). 
\end{eqnarray*}
%%%%%%%%%%%%%
A similar result holds for $x\in \widehat{\Sigma}^{(g_2)}$, $y\in \widehat{\Sigma}^{(g_{1} )}$. 
\hfill $\square $
\medskip

%%%%%%%%%%%%%%%%%%%%%%%%%%%%%%%%%%%%%%%%%%%%%%%%%%%%%%%%%%%%%%%%%%%%%%%%%%%%%%%%%%%
We next compute the explicit form of the moment matrix $X_{ab}$ in terms of
the moments $F_a$ of $S^{(g_{a})}(x,y)$. 
%%%%%%%%%%%%%%%%%%%%%%%%%%%%%%%%%%%%%%%%%%%%%%%%%%%%%%
%%
It is useful to introduce infinite block matrices 
%%%%%%%%%%%%%%%%%%%%%%%%%%%%%%%%%%%%
\begin{eqnarray}
&&X=\left( 
\begin{array}{cc}
X_{11} & X_{12} \\ 
X_{21} & X_{22}
\end{array}
\right) ,
\qquad 
F=\left( 
\begin{array}{cc}
F_{1} & 0 \\ 
0  & F_{2}
\end{array}
\right),
\nonumber\\
&&\Xi=\left( 
\begin{array}{cc}
0 & \xi I \\ 
-\xi I & 0
\end{array}
\right),
\qquad 
Q=F\Xi =\left( 
\begin{array}{cc}
0 & \xi F_{1} \\ 
-\xi F_{2} & 0
\end{array}
\right). 
  \label{XFQ_def}
\end{eqnarray}
Then one finds:
%%%
\begin{proposition}
\label{Prop_F1F2}
$X$ is given by
\begin{equation}
X=(I-Q)^{-1}F,
\label{XQsol}
\end{equation} 
where $(I-Q)^{-1}=\sum\limits_{n\geq 0}Q^n$
is convergent for $|\epsilon|<\vert r_{1}  r_2\vert$. 
\end{proposition}
%%%%%%%%%%%%%%%%%%%%%%%%%%%%%%%%%%%%%%%%%%%%%%%%%%%%%%%%%%%%%

\noindent \textbf{Proof.}  
%%%%%%%%%%%%%%%%%%%%%%%%%%%%%%%%%%%%%%%%%%%%%%%%%%%%%%%%%%%%%%%%%
Using (\ref{Prop_epsinteq1}) we find $X_{11}(k,l)-F_{1} (k,l)$ is given by
%%%%%%%%%%%%%%%%%%%%%%%%%%%%%%%%%
\begin{eqnarray*}
&&-\frac{\epsilon^{\half(k+l-1)}}{(\tpi )^{3}}
\oint\limits_{\mathcal{C}_{1}(x)}
\oint\limits_{\mathcal{C}_{1}(z_{1})}
S^{(g_{1})}(x,z_{1} )x^{-k}dx^{\half}
\oint\limits_{\mathcal{C}_{1}(y)}
S^{(g_{1}+g_2)}(z_{1} ,y)y^{-l}dy^{\half}
\\
&&=-\frac{\epsilon^{\half(k+l-1)}}{(\tpi )^{3}}
\sum_{m\ge 0}
\left(
 \epsilon^{-\frac{m}{2} +\frac{1}{4}} 
 \oint\limits_{ \mathcal{C}_{1}(x)} h_{1}(m,x) x^{-k}dx^{\half}\right.
\\
&&\left. \cdot
\oint\limits_{\mathcal{C}_{1}(z_{1})}
\oint\limits_{\mathcal{C}_{1}(y)}
S^{(g_{1}+g_2)}(z_{1} ,y)z_{1} ^{m-1}y^{-l}dy^{\half}
\right)
\\ 
&&=\xi\sum_{m\ge 0} 
\left(
\frac{\epsilon^{\frac{k}{2} -\frac{1}{4}}}{\tpi } 
\oint\limits_{ \mathcal{C}_{1}(x)} h_{1}(m,x) x^{-k}dx^{\half}
\right.
\\
&&\left.
\cdot
\frac{\epsilon^{\half(m+l-1)}}{(\tpi )^{2}}\oint\limits_{\mathcal{C}_{2}(z_{2})}
\oint\limits_{\mathcal{C}_{1}(y)}
S^{(g_{1}+g_2)}(z_2,y)z_2^{-m}y^{-l}dy^{\half}
\right)
\\  
&&=\xi \left(F_{1}  X_{21}\right)(k,l),
\end{eqnarray*}
%%%%%%%%%%%%%%%%%%%%%%%%%%%%%%%%%%%%%%%%%%%%%%%%%%%%%%%%%%%%%%
using (\ref{htoS}) and (\ref{dz1dz2}). Similarly we find $X_{22}=F_2-\xi F_2 X_{12}$ so that
$X_{aa}=(F+QX)_{aa}$ using (\ref{XFQ_def}). 
%%%%%%%%%%%%%%%%%%%%%%%%%%%%%%%%%%%%%%%%%%%%%%%%%%%%%%%%%%%%%%%%%%%%%%%%%%%%%
A similar calculation of $X_{12}$ and $X_{21}$ leads to $X_{a\bar{a}}=(QX)_{a\bar{a}}$.
%%%%%%%%%%%%%%%%%%%%%%%
%% 
These combine to give  
$(I-Q)X=F$ which implies (\ref{XQsol}) provided $(I-Q)^{-1}=\sum_{n\ge 0}Q^n$ converges.  	
But (\ref{XQsol}) can be rewritten 
\begin{equation}
X=\sum_{n\ge 1}Q^n\Xi.
\label{XabQ}
\end{equation}
%%
%%%%%%%%%%%%%% 
By Theorem \ref{theoremomg1g2holo}, $X_{ab}(k,l)$ has a convergent series expansion in $\epsilon^{\half}$ for $|\epsilon |<r_{1}r_{2}$. But $\sum_{n=1}^NQ^n=O(\epsilon^{\half N})$ so that (\ref{XabQ}) holds to all orders in $\epsilon^{\half}$. Hence $(I-Q)^{-1}$ converges for $|\epsilon |<r_{1}r_{2}$ and the proposition holds. 
\hfill $\square$ 

Propositions \ref{Prop_hXh} and \ref{Prop_F1F2} imply
\begin{theorem} 
\label{Theorem_Shh}
$S^{(g_{1}+g_{2})}(x,y)$ is given by 
\begin{equation*}
S^{(g_{1}+g_{2})}(x,y)
=\delta_{ab} S^{(g_{a})}(x,y)
 +h_{a}(x) 
 \left(\Xi(I-Q)^{-1}\right)_{ab}  
 \bar{h}_{b}^{T}(y),
\end{equation*}
for $x \in \widehat{\Sigma}^{(g_a)}, y \in \widehat{\Sigma}^{(g_b)}$. Equivalently, 
\begin{equation*}
S^{(g_{1}+g_{2})}(x,y)
=\left\{ 
\begin{array}{l}
S^{(g_{a})}(x,y)
 + h_{a}(x) 
 \left(I-F_{\bar{a}}F_a\right)^{-1}F_{\bar{a}}  
 \bar{h}_{a}^{T}(y),  \quad  \mbox{ for } x,y \in 
\widehat{\Sigma}^{(g_a)}, 
\\
  \xi (-1)^{\bar{a}} h_{a}(x)
 \left(  I- F_{\bar{a}}F_a\right)^{-1} 
 \bar{h}_{\bar{a}}^{T}(y), \quad 
\mbox{ for } x\in \widehat{\Sigma}^{(g_a)},
\; y\in \widehat{\Sigma}^{(g_{\bar{a}})}. \quad \square
\end{array} 
\right. 
\end{equation*}
\end{theorem} 

\begin{remark}\label{rem_epsilon_root_branch}
Note that $S^{(g_{1}+g_{2})}(x,y)$ is even (odd) in $\epsilon^{\half}$ for 
$x,y \in \widehat{\Sigma}^{(g_a)}$ (respectively,  for $x\in \widehat{\Sigma}^{(g_a)}$, $y\in \widehat{\Sigma}^{(g_{\bar{a}})}$). Thus $S^{(g_{1}+g_{2})}(x,y)$ is invariant under a Dehn twist $\epsilon\rightarrow e^{\tpi }\epsilon$ with  $\xi\rightarrow -\xi$ from (\ref{dz1dz2}).

\end{remark} 
\medskip
Similarly to ref. \cite{MT1} we define the determinant of $I-Q$ 
as a formal power series in $\epsilon^{\half}$ by
\begin{equation*}
\log \det\left(I - Q\right) = \mathrm{Tr} \log\left(I - Q\right)=-\sum_{n\ge 1}\frac{1}{n} \mathrm{Tr}(Q^n).
\end{equation*}
Clearly $\mathrm{Tr}(Q^{2k})=2\mathrm{Tr}\left((F_{1} F_2)^k\right)$ for $k\ge 0$ whereas $\mathrm{Tr}(Q^n)=0$ for $n$ odd. Furthermore, from (\ref{Fdef}) the diagonal terms $(F_{1} F_2)^k$ have integral power series in $\epsilon$. Thus it follows that
\begin{lemma}\label{Lemma_logdet}
$\det\left(I - Q\right)=\det\left(I - F_{1} F_2\right)$ and is a formal power series in $\epsilon$.
\end{lemma}  
The determinant has the following holomorphic properties:
%%%%%%%%%%%%%%%%%%%%%%%%%%%%%%%%%%%%%%%%%%%
\begin{theorem}
\label{theorem_Det} 
$\det \left(I-Q\right)$ is non-vanishing and holomorphic
in $\epsilon $ for $|\epsilon |<r_{1}r_{2}$.
\end{theorem}

\noindent \textbf{Proof.} The proof follows a similar argument to Theorem 2 of ref. \cite{MT1}.
Let $S^{(g_{1}+g_{2})}(z_{1},z_{2})=f(z_{1},z_{2},\epsilon )dz_{1}^{\half}dz_{2}^{\half}$ for $|z_{a}|\leq
r_{a} $ where $f(z_{1},z_{2},\epsilon )$ is holomorphic in $\epsilon^{\half} $ for $%
|\epsilon |\leq r$ for $r<r_{1}r_{2}$ from Theorem \ref{theoremomg1g2holo}.
Apply Cauchy's inequality to the coefficients of
$f(z_{1},z_{2},\epsilon )
=\sum_{n\geq 0}f_{n}(z_{1},z_{2})\epsilon^{{\frac{n}{2}}}$ 
to find 
\begin{equation}
|f_{n}(z_{1},z_{2})|\leq Mr^{-{\frac{n}{2}}},  \label{Cauchy}
\end{equation}%
for $M=\sup_{|z_{a}|\leq r_{a},|\epsilon |\leq r}|f(z_{1},z_{2},\epsilon )|$. 
Consider 
\begin{equation}
\mathcal{I}=\frac{1}{(\tpi )^{2}}\oint_{\mathcal{C}_{r_{1}}(z_{_{1}})}%
\oint_{\mathcal{C}_{r_{2}}(z_{_{2}})}
S^{(g_{1}+g_{2})}(z_{1},z_{2})
\left(1-\frac{\epsilon }{z_{1}z_{2}}\right)^{-1}
dz_{1}^{\half}dz_{2}^{\half},
\label{Siberian}
\end{equation}%
for $\mathcal{C}_{r_{a}}(z_{a})$ the contour with $|z_{a}|=r_{a}$. Then
using (\ref{Cauchy}) we find 
\begin{eqnarray*}
|\mathcal{I}| \leq M.\sum_{n\geq 0}
\left(\frac{|\epsilon |}{r}\right)^{{\frac{n}{2}}}.
\left|1-\frac{|\epsilon |}{r_{1}r_{2}}\right|^{-1}.r_{1}r_{2},
\end{eqnarray*}%
i.e. $\mathcal{I}$ is absolutely convergent and thus holomorphic in $%
\epsilon^{{\half}} $ for $|\epsilon |<r<r_{1}r_{2}$. 
Since $|z_{1}z_{2}|=r_{1}r_{2}$ we may alternatively expand in $\epsilon/z_{1}z_{2}$ to obtain 
\begin{eqnarray*}
\mathcal{I} &=&\sum_{k\geq 1}\epsilon^{k}\frac{1}{(\tpi )^{2}}%
\oint_{\mathcal{C}_{r_{1}}(z_{_{1}})}\oint_{\mathcal{C}_{r_{2}}(z_{_{2}})}%
S^{(g_{1}+g_{2})}(z_{1},z_{2})z_{1}^{-k}z_{2}^{-k}
dz_{1}^{\half}dz_{2}^{\half} \\
&=&\epsilon^{{\half}}\mathrm{Tr} X_{12},
\end{eqnarray*}%
where $\mathrm{Tr}X_{12}=\sum_{k\geq 1}X_{12}(k,k)$.
But \eqref{XQsol} implies
\begin{equation*}
\mathrm{Tr}X_{12}=\xi\sum_{n\geq 1}\mathrm{Tr}((F_{1}F_{2})^{n}),
\end{equation*}%
which is absolutely convergent for $|\epsilon |<r_{1}r_{2}$. Hence we find 
\begin{equation*}
\mathrm{Tr}\log (I-F_{1}F_{2})=-\sum_{n\geq 1}\frac{1}{n}\mathrm{Tr}((F_{1}F_{2})^{n}),
\end{equation*}%
is also absolutely convergent for $|\epsilon |<r_{1}r_{2}$. Thus $\det (I-Q)=\det (I-F_{1}F_{2})$ is non-vanishing and holomorphic for $|\epsilon |<r_{1}r_{2}$. \hfill $\square$

%%%%%%%%%%%%%%%%%%%%%%%%%%%%%%%%%%%%%%%%%%%%%%%%%%%%%%%%%%%%%%%%%%%%%%%
%%%%%%%%%%%%%%%%%%%%%%%%%%%%%%%%%%%%%%%%%%%%%%%%%%%%%%%%%%%%%%%%%%%%%%%
\subsection{Sewing Two Tori}
\label{subsect_genus_two_eps}
%%
%%%%%%%%%%%%%%%%%%%%%%%%%%%%%%%%%%%%%%%%%%%%%%%%%%%%%%%%%%%%%%%%%%%%%%%
Consider the genus two surface  formed by sewing two oriented tori  $\Sigma_a^{(1)}=\mathbb{C}/{\Lambda}_{a}$ for $a=1,2$, 
and lattice 
${\Lambda}_{a}=\tpi (\mathbb{Z}\tau _{a}\oplus \mathbb{Z})$ for $\tau _{a}\in 
\mathbb{H}_{1}$. This is discussed at length in \cite{MT1}. 
 For local coordinate $z_{a}\in \mathbb{C}/{\Lambda}_{a}$
consider the closed disk $\left\vert z_{a}\right\vert \leq r_{a}$ which is
contained in $\Sigma^{(1)}_a$ provided\ $r_{a}<\half D(q_{a})$ where 
\begin{equation*}
D(q_{a})=\min_{\lambda \in {\Lambda}_{a}, \lambda \neq 0}|\lambda |,
\end{equation*}
is the minimal lattice distance. From Subsection \ref{Subsec_eps_sew} we obtain a genus two Riemann surface  $\Sigma^{(2)}$ parameterized 
by the domain  
\begin{equation}
\mathcal{D}^{\epsilon }=\{(\tau _{1},\tau _{2},\epsilon )\in \mathbb{H}_{1}
\mathbb{\times H}_{1}\mathbb{\times C}\ :\ |\epsilon |<\frac{1}{4}
D(q_{1})D(q_{2})\}.  
\label{Deps}
\end{equation}
$\mathcal{D}^{\epsilon }$ is preserved under the action of 
$G\simeq (SL(2,\mathbb{Z})$ $\times SL(2,\mathbb{Z}))\rtimes \mathbb{Z}_{2}$, the direct
product of the left and right torus modular groups, which are interchanged
upon conjugation by an involution $\beta $ as follows
\begin{eqnarray}
\gamma _{1}(\tau _{1},\tau _{2},\epsilon ) &=&\left(\frac{a_{1}\tau _{1}+b_{1}}{
c_{1}\tau _{1}+d_{1}},\tau _{2},\frac{\epsilon }{c_{1}\tau _{1}+d_{1}}\right), 
\notag \\
\gamma _{2}(\tau _{1},\tau _{2},\epsilon ) &=&\left(\tau _{1},\frac{a_{2}\tau
_{2}+b_{2}}{c_{2}\tau _{2}+d_{2}},\frac{\epsilon }{c_{2}\tau _{2}+d_{2}}\right), 
\notag \\
\beta (\tau _{1},\tau _{2},\epsilon ) &=&(\tau _{2},\tau _{1},\epsilon ),
\label{eq: GDeps}
\end{eqnarray}%
for $(\gamma _{1},\gamma _{2})\in SL(2,\mathbb{Z})\times SL(2,\mathbb{Z})$
with $\gamma _{i}=\left( 
\begin{array}{cc}
a_{i} & b_{i} \\ 
c_{i} & d_{i}%
\end{array}%
\right) $.
\medskip

There is a natural injection $G\rightarrow Sp(4,\mathbb{Z})$ in which the
two $SL(2,\mathbb{Z})$ subgroups are mapped to 
\begin{equation}
\Gamma _{1}=\left\{ \left[ 
\begin{array}{cccc}
a_{1} & 0 & b_{1} & 0 \\ 
0 & 1 & 0 & 0 \\ 
c_{1} & 0 & d_{1} & 0 \\ 
0 & 0 & 0 & 1%
\end{array}%
\right] \right\} ,\;\Gamma _{2}=\left\{ \left[ 
\begin{array}{cccc}
1 & 0 & 0 & 0 \\ 
0 & a_{2} & 0 & b_{2} \\ 
0 & 0 & 1 & 0 \\ 
0 & c_{2} & 0 & d_{2}%
\end{array}%
\right] \right\},  \label{eq: G1G2}
\end{equation}%
and the involution is mapped to 
\begin{equation}
\beta =\left[ 
\begin{array}{cccc}
0 & 1 & 0 & 0 \\ 
1 & 0 & 0 & 0 \\ 
0 & 0 & 0 & 1 \\ 
0 & 0 & 1 & 0%
\end{array}%
\right].  \label{eq: beta}
\end{equation}%
$G$ also has a natural action on $\mathbb{H}_{2}$ as given in (\ref{eq:modOmega}) which is compatible with respect to $\Omega^{(2)}$ as a  function of $(\tau _{1},\tau _{2},\epsilon )$  \cite{MT1}.

%%%%%%%%%%%%%%%%%%%%%%%%%%%%%%%%%%%%%%%%%%%%%%%%%%%%%%%%%
\medskip
The Szeg\"o kernel on the torus ${\Sigma}_a^{(1)}$ is given by   
\begin{equation*}
S^{(1)} \left[ { \theta}_a  \atop { \phi }_a \right]( x,y|\tau_a)   
=P_{1}  \left[ { \theta}_a  \atop { \phi }_a \right]( x-y, \tau_a) dx^{\half} dy^{\half},
\end{equation*}
 from (\ref{S1}). 
It is straightforward to compute the moment matrix $F_a$  of (\ref{Fdef}). Using the Laurant expansion (\ref{P1zn}) we find  
\begin{equation}
P_{1}\left[{\theta} \atop  {\phi }\right] 
(x-y,\tau )=\frac{1}{x-y}+\sum_{k,l\geq 1}
C\left[{\theta} \atop  {\phi }\right]
(k,l)x^{k-1}y^{l-1},\label{P1xy}
\end{equation}%
where for $k,l\geq 1$ we define 
\begin{equation}
C\left[{\theta} \atop  {\phi }\right] 
 (k,l,\tau )=(-1)^{l}\binom{k+l-2}{k-1}
 E_{k+l-1}\left[{\theta} \atop  {\phi }\right] 
 (\tau ),\label{Cdef}
\end{equation}
for twisted  Eisenstein series \eqref{Endef}. Then it follows that
\begin{equation}
F_{a}\left[ { \theta}_a  \atop { \phi }_a \right] ( k,l,\tau _{a},\epsilon )
=\epsilon^{\half(k+l-1)}
 C\left[ { \theta}_a  \atop { \phi }_a \right]  (  k,l,\tau _{a}). 
  \label{FaCdef}
\end{equation} 
We also have the analytic expansion
\begin{eqnarray}
P_{1}  \left[ { \theta}  \atop { \phi } \right]( x-y, \tau)
&=&\sum_{k\ge 0} P_k\left[ { \theta}  \atop { \phi } \right]( x, \tau) y^{k-1},\label{P1xyk}
\end{eqnarray}
for $P_{k}\left[ { \theta}  \atop { \phi }\right]( z, \tau)
=\frac{(-1)^{k-1}}{(k-1)!}\partial^{k-1}_z P_{1} \left[ { \theta} \atop { \phi } \right]( z, \tau)$. Then we find 
\begin{equation}
h_a\left[ { \theta_a}  \atop { \phi_a }\right](k,x, \tau_a,\epsilon)=
\epsilon^{\frac{k}{2}-\frac{1}{4}} 
P_{k}\left[ { \theta_a}  \atop { \phi_a }\right]
(x, \tau_a)dx^{\half}.
\label{haCdef}
\end{equation}
Using these results we may therefore determine the explicit form for $S^{(2)}
\left[
{\theta^{(2)}} \atop {\phi^{(2)}}	
\right]$ on $\mathcal{D}^{\epsilon }$ via Theorem \ref{Theorem_Shh}.  
\medskip

One may also confirm that $S^{(2)}$ satisfies  the modular invariance property of (\ref{Szmod})  under the group $G$ generated by $\gamma_i,\beta$ of (\ref{eq: G1G2}) and (\ref{eq: beta}) with  
\begin{equation}
S^{(2)}\left(\gamma
\left[
{\theta^{(2)}} \atop {\phi^{(2)}}	
\right]\right)
(\gamma x, \gamma y \vert \gamma (\tau _{1},\tau _{2},\epsilon ))
=
S^{(2)}
\left[
{\theta^{(2)}} \atop {\phi^{(2)}}	
\right]
(x,y\vert\tau _{1},\tau _{2},\epsilon ),
\label{eq:S2mod}
\end{equation}
where
\begin{equation*}
\gamma_{1} 
\left[
\begin{array}{c}
\theta_{1} \\
\theta_2\\
\phi_{1} \\
\phi_2
\end{array}	
\right]=
\left[
\begin{array}{c}
\theta_{1} ^{a_{1} }\phi_{1} ^{b_{1} }\\
\theta_2\\
\theta_{1} ^{c_{1} }\phi_{1} ^{d_{1} }\\
\phi_2
\end{array}	
\right],\quad
\gamma_2\left[
\begin{array}{c}
\theta_{1} \\
\theta_2\\
\phi_{1} \\
\phi_2
\end{array}	
\right]
=
\left[
\begin{array}{c}
\theta_{1} \\
\theta_2^{a_2}\phi_2^{b_2}\\
\phi_{1} \\
\theta_2^{c_2}\phi_2^{d_2}\\
\end{array}	
\right],\quad
\beta
\left[
\begin{array}{c}
\theta_{1} \\
\theta_2\\
\phi_{1} \\
\phi_2
\end{array}	
\right]
=
\left[
\begin{array}{c}
\theta_2\\
\theta_{1} \\
\phi_2\\
\phi_{1} 
\end{array}	
\right],
\end{equation*} 
and 
\begin{eqnarray*}
\gamma _{a}x &=& \left\{
\begin{array}{cl}
\frac{x}{c_{a}\tau _{a}+d_{a}},	& \mbox{for } x\in \widehat{\Sigma}^{(1)}_a,
\\
x, & \mbox{for } x\in \widehat{\Sigma}^{(1)}_{\bar{a}},
\end{array}
\right. 
\end{eqnarray*}
and where for $x=\tpi   \left(u + v\tau_a \right)\in \widehat{\Sigma}^{(1)}_a$ with $0\le u,v <1$ we  define  
$\beta  x = \tpi   \left(u + v\tau_{\bar{a}} \right)$. Finally, we note that $\det \left(I-Q\right)=\det \left(I-
F_1\left[ { \theta}_1  \atop { \phi }_1 \right] 
F_2\left[ { \theta}_2  \atop { \phi }_2 \right]\right)$ is also $G$ invariant.
  
\section{The Szeg\"o kernel on a Self-Sewn Riemann Surface}
\label{Sect_Rho_g}
%%%%%%%%%%%%%%%%%%%%%%%%%%%%%%%%%%%%%%%%%%%%%%%%%%%%%%%%%

\subsection{The $\rho$-Formalism Sewing Scheme}
We now consider the construction of the Szeg\"o kernel on a Riemann surface $\Sigma^{(g+1)}$ formed by self-sewing a handle to a Riemann surface $\Sigma^{(g)}$ of genus $g$. We begin by reviewing the Yamada formalism \cite{Y} in this scheme which, following \cite{MT1}, we refer to as the $\rho$-formalism. 
Consider a Riemann surface $\Sigma^{(g)}$ of genus $g$ and let $z_{1},z_{2}$ be local coordinates in the neighborhood of two separated points $p_{1}$ and $p_{2}$.
Consider two disks $\left\vert z_{a}\right\vert
\leq r_{a}$, for $r_{a}>0$ and $a=1,2$. Note that 
$r_{1},r_{2}$ must be sufficiently small to ensure that the disks do not
intersect. Introduce a complex parameter ${\rho }$ where $|{\rho }|\leq
r_{1}r_{2}$ and excise the disks
\begin{equation*}
\{z_{a}:\, \left\vert z_{a}\right\vert <|\rho |r_{\bar{a}}^{-1}\}\subset 
\Sigma^{(g)}, 
\end{equation*} 
to form a twice-punctured surface 
\begin{equation*}
\widehat{\Sigma}^{(g)}=\Sigma^{(g)}\backslash \bigcup_{a=1,2}\{z_{a}:\, \left\vert
z_{a}\right\vert <|\rho |r_{\bar{a}}^{-1}\}.
\end{equation*}%
As before, we use the convention 
%%(\ref{bardef})
$\bar 1=2$, $\bar 2=1$.
 We define annular regions $%
\mathcal{A}_{a}\subset \widehat{\Sigma}^{(g)}$ with $\mathcal{A}_{a}=
\{z_{a}:\, |{
\rho }|r_{\bar{a}}^{-1}\leq \left\vert z_{a}\right\vert \leq r_{a}\}$ and
identify them as a single region $\mathcal{A}=\mathcal{A}_{1}\simeq \mathcal{%
A}_{2}$ via the sewing relation 
\begin{equation}
z_{1}z_{2}=\rho ,  \label{rhosew}
\end{equation}%
to form a compact Riemann surface 
$\Sigma^{(g+1)}=\widehat{\Sigma}^{(g)}\backslash \{\mathcal{A}%
_{1}\cup \mathcal{A}_{2}\}\cup \mathcal{A}
$ of genus $g+1$. The sewing
relation (\ref{rhosew}) can be considered to be a parameterization of a cylinder connecting the punctured Riemann surface to itself. 

%%%%%%%%%%%%%%%%%%%%%%%%%%%%%%%%%%%%%%%%%%%%%%%%%%%%%%%%%%%%%%%%%%%%%%%%%%%%%%%%%%%%%%%%%%%%
%%
In the $\rho$-formalism we define a standard basis of cycles $\{a_{1},b_{1},\ldots, a_{g+1}, b_{g+1}\}$ on $\Sigma^{(g+1)}$ where the set 
$\{a_{1}, b_{1}, \ldots, a_{g}, b_{g}\}$ is the original basis on $\Sigma^{(g)}$. Let $\mathcal{C}_{a}(z_{a})\subset \mathcal{A}_{a}$ denote a closed anti-clockwise contour parameterized by $z_{a}$ surrounding the puncture at $z_{a}=0$. Clearly $\mathcal{C}_{2}(z_{2})\sim -\mathcal{C}_{1}(z_{1})$ on applying the sewing relation \eqref{rhosew}. We then define the cycle $a_{g+1}$ to be $\mathcal{C}_{2}(z_{2})$ and define the cycle $b_{g+1}$ to be a path chosen in $\widehat{\Sigma}^{(g)}$ between identified points $z_{1}=z_{0}$ and $z_{2}=\rho /z_{0}$ on the sewn surface. 

As in the $\epsilon$-formalism, the normalized differential of the second kind $\omega^{(g+1)}$, the holomorphic 1-forms $\nu_i^{(g+1)}$ and the period matrix 
$\Omega^{(g+1)}$ can be computed in terms of data coming from $\Sigma^{(g)}$ \cite{Y},\cite{MT1} to find
%%%%%%%%%%%%%%%%%%%%%%%%%%%%%%%%%%%%%%%%%%%%%%%%%
\begin{theorem}
\label{theoremomgplus} 
$\omega^{(g+1)}$, $\nu_{i}^{(g+1)}$ and $\Omega_{ij}^{(g+1)}$ for $(i,j)\ne (g+1,g+1)$ are holomorphic in $\rho $ for 
$|\rho |<r_{1}r_{2}$ with
\begin{eqnarray}
&&\omega^{(g+1)}(x,y)=
\omega^{(g)}(x,y)+O(\rho),
\label{omega_rho}
\\
&&\nu^{(g+1)}_{i}(x)=
\nu^{(g)}_{i}(x)+O(\rho),\quad i=1, \ldots,  g
\label{nui_rho}
\\
&&\nu^{(g+1)}_{g+1}(x)=
\omega^{(g)}_{p_2-p_{1} }(x)+O(\rho),
\label{nugp_rho}
\\
&&\Omega^{(g+1)}_{ij}=
\Omega^{(g)}_{ij}+O(\rho),\quad i,j=1, \ldots, g
\label{Omij_rho}
\\
&&\Omega^{(g+1)}_{i,g+1}=
\frac{1}{\tpi }\int_{p_{1} }^{p_2}\nu^{(g)}_{i}+O(\rho),\quad i=1, \ldots, g,
\label{Omigp_rho}
\end{eqnarray}
for $x,y\in \widehat{\Sigma}^{(g)}$. 
$e^{\tpi  \Omega^{(g+1)}_{g+1,g+1}}$ is holomorphic in $\rho $ for 
$|\rho |<r_{1}r_{2}$ with
\begin{eqnarray}
e^{\tpi  \Omega^{(g+1)}_{g+1,g+1}}=-
\frac{\rho}{K_0^2}
\left( 1+O(\rho)\right),
\label{Omgpgp_rho}
\end{eqnarray}  
where $K_0=K^{(g)}(z_{1} =0,z_2=0)$ for  $E^{(g)}(z_{1} ,z_2)=K^{(g)}(z_{1} ,z_2)dz_{1} ^{-\half}dz_2^{-\half}$ expressed in terms of the local coordinates $z_{1} ,z_2$. \hfill $\square$  
\end{theorem} 
%%%%%%%%%%%%%%%%%%%%%%%%%%%%%%%%%%%%%%%%%%%%%%%%%%%%%%%%%%%%%%%%%%%%%
\subsection{Szeg\"o Kernel in the $\rho$-Formalism} 
%%%%%%%%%%%%%%%%%%%%%%%%%%%%%%%%%%%%%%%%%%%%%%%%%%%%%%%%%%%%%%%%%%%%%
\label{relationslambda1}
We now determine the  Szeg\"o kernel $S^{(g+1)}(x,y)=S^{(g+1)}\left[ {\theta}^{(g+1)}  \atop {\phi}^{(g+1)} \right](x,y)$ on the sewn Riemann surface $\Sigma^{(g+1)}$ in terms of genus $g$ data
together with the multiplier parameters associated with the handle cycles.
The $S^{(g+1)}$ multipliers \eqref{eq:periodicities} on the cycles $a_i,b_i$ for $i=1,\ldots, g$ are determined by the multipliers of $S^{(g)}$ with
$\phi_i^{(g+1)}=\phi_i^{(g)}$ and $\theta_i^{(g+1)}=\theta_i^{(g)}$ 
i.e. $\alpha_i^{(g+1)}=\alpha_i^{(g)}$ and $\beta_i^{(g+1)}=\beta_i^{(g)}$.  
The remaining two multipliers associated with the cycles $a_{g+1}$ and $b_{g+1}$ 
\begin{eqnarray}
\phi_{g+1}&=&\phi_{g+1}^{(g+1)}=-e^{\tpi \alpha_{g+1}^{(g+1)}},\label{phigp}\\
\theta_{g+1}&=&\theta_{g+1}^{(g+1)}=-e^{-\tpi \beta_{g+1}^{(g+1)}},\label{thetagp}
\end{eqnarray}
must be  additionally specified so that
\begin{eqnarray}
 S^{(g+1)}(e^{\tpi }x_a,y) &=& -\phi_{g+1}^{ a-\bar{a}} \; 
 S^{(g+1)}(x_{a},y),    \label{agplus_mult}\\
 S^{(g+1)}(x_a,y) &=& -\theta_{g+1}^{ a-\bar{a}} \; 
 S^{(g+1)}(x_{\bar{ a}},y),    \label{bgplus_mult}
\end{eqnarray}
for $x_{a}\in \mathcal{A}_{a}$ and $x_{\bar{a}}\in\mathcal{A}_{\bar a}$.

We next consider the analogue of Theorem \ref{theoremomg1g2holo} concerning the holomorphicity of $S^{(g+1)}$ as a function of $\rho$. It is convenient to define $\kappa\in \left[-\half,\half\right)$ by
$\phi_{g+1}=-e^{\tpi \kappa}$ i.e. $\kappa=\alpha_{g+1}^{(g+1)} \mod 1$. 
We then find  
\begin{theorem}
\label{theoremomgp1holo11}
 $S^{(g+1)}$ is holomorphic in $\rho $ for $|\rho|<r_{1}r_{2}$ with 
\begin{equation}
S^{(g+1)}(x,y)=S_{\kappa}^{(g)}(x,y)+O(\rho),
\label{eq:SgpSkappa}
\end{equation} 
for $x,y\in \widehat{\Sigma}^{(g)}$ where $S_{\kappa}^{(g)}(x,y)$ is defined as follows: For $\kappa\neq -\half$ 
\begin{eqnarray}
S_{\kappa}^{(g)}(x,y)&=&
\frac{U(x,y)^{\kappa} \vartheta \left[ {{\alpha}^{(g)} \atop {\beta^{(g)}}} \right] 
\left( \int_{y}^{x}\nu^{(g)}  +\kappa z_{p_{1} ,p_2}\vert \Omega^{(g)}\right)}
{E^{(g)}(x, y)\vartheta \left[ {{\alpha^{(g)}} \atop {\beta^{(g)}}}\right] 
\left(\kappa z_{p_{1} ,p_2}\vert \Omega^{(g)}\right) },
\label{eq:Skappa}
\end{eqnarray}
where
\begin{equation}
U(x,y)=\frac { E^{(g)}(x, p_2) E^{(g)}(y, p_{1} )  } 
{ E^{(g)} (x, p_{1} )  { E^{(g)}(y, p_2) } }, 
\label{eq:U}
\end{equation}
for prime form $E^{(g)}$ and where
\begin{equation}
z_{p_{1} ,p_2}=\int_{p_{1} }^{p_2} \nu^{(g)},
\label{eq:z12}
\end{equation}
 for holomorphic 1-forms $\nu^{(g)}$. For $\kappa=-\half$ then $S_{-\half}^{(g)}(x,y)$ is given by
\begin{eqnarray}
&
\left(
\frac{U(x,y)^{\half}}{E^{(g)}(x,y)}
\vartheta \left[ {{\alpha}^{(g)} \atop {\beta^{(g)}}} \right] 
\left( \int_{y}^{x}\nu^{(g)}  +\half z_{p_1,p_2}\vert \Omega^{(g)}\right)
\right.&\notag \\
&\left.
-\theta_{g+1}\frac{U(x,y)^{-\half}}{E^{(g)}(x,y)}
\vartheta \left[ {{\alpha}^{(g)} \atop {\beta^{(g)}}} \right] 
\left( \int_{y}^{x} \nu^{(g)} -\half z_{p_1,p_2}\vert \Omega^{(g)}\right)
\right).&\notag\\
&
\left(
\vartheta \left[ {{\alpha^{(g)}} \atop {\beta^{(g)}}} \right]\left(\half z_{p_1,p_2}\vert \Omega^{(g)}\right) 
-\theta_{g+1}\vartheta \left[ {{\alpha^{(g)}} \atop {\beta^{(g)}}} \right]\left(-\half z_{p_1,p_2}\vert \Omega^{(g)}\right)
\right)^{-1} .&
\label{eq:Shalf} 
\end{eqnarray}
\end{theorem}
\noindent{\bf Proof.}  
We firstly note that from \eqref{eq:Prime_omega} it follows that 
\begin{equation}
E^{(g+1)}(x,y)=E^{(g)}(x,y)+O(\rho).\label{E_rhoplus}
\end{equation} 
From Theorem \ref{theoremomgplus} we may expand the genus $g+1$ theta series to leading order in $\rho$ for $|\rho|<r_{1}r_{2}$ as follows
\begin{eqnarray*}
&&\vartheta\left[{ {\alpha^{(g+1)}} \atop {\beta^{(g+1)}} }\right]
\left( \int_{y}^{x} {\nu}^{(g+1)} |\Omega^{(g+1)}\right)
=
\sum_{ m \in {\mathbb{Z}}^{g} }\sum_{ n \in {\mathbb{Z}}}
\left( -
\frac{\rho}{K_{0}^2} \right)^{\half(n+\alpha_{g+1}^{(g+1)})^2}.  \\
&&\exp\left(  i \pi (m + \alpha^{(g)}).\Omega^{(g)} .(m+{\alpha^{(g)}}) + 
(m+ {\alpha^{(g)}}). ({\int_{y}^{x} {\nu}^{(g)}+ \tpi  {\beta^{(g)}}})+\right. \\
&&\left. (n+{\alpha_{g+1}^{(g+1)}})
\left[(m + \alpha^{(g)}).\int_{p_{1} }^{p_2} {\nu}^{(g)} 
+{\int_{y}^{x} {\omega}_{p_2-p_{1} }^{(g)}+ \tpi  {\beta_{g+1}^{(g+1)}}}
\right]
\right).\\
&&(1+O(\rho) ),
\end{eqnarray*} 
Clearly
$|n+\alpha_{g+1}^{(g+1)}|\ge |\kappa|$. For $\kappa\neq-\half$ it follows that this lower bound is satisfied for one value of $n$ so that
\begin{eqnarray*}
&&\vartheta\left[{ {\alpha^{(g+1)}} \atop {\beta^{(g+1)}} }\right]
\left( \int_{y}^{x} {\nu}^{(g+1)} |\Omega^{(g+1)}\right)
=
\left( -\frac{\rho}{K_0^2} \right)^{\half\kappa^2}(-\theta_{g+1})^{-\kappa}U(x,y)^{\kappa}
\\
&&
\sum_{ m \in {\mathbb{Z}}^{g} }
\exp\left( i \pi (m + \alpha^{(g)}).\Omega^{(g)} .(m+{\alpha^{(g)}}) \right.\\ 
&&+\left. (m+ {\alpha^{(g)}}). ({\int_{y}^{x} {\nu}^{(g)}+\kappa z_{p_1,p_2}+ \tpi  {\beta^{(g)}}})
\right).(1+O(\rho)),
\end{eqnarray*}
for $z_{p_1,p_2}$ of \eqref{eq:z12} and  
where from \eqref{omprime} 
\begin{equation*}
\int\limits_{y}^{x} \omega^{(g)}_{p_2-p_{1} }=  
\int\limits_{y}^{x} \int\limits_{p_{1} }^{p_2} \omega^{(g)}(\cdot,\cdot)
=
\log  U(x,y),
%%%%    
%%
\end{equation*}
for $U(x,y)$ of \eqref{eq:U}. Therefore
\begin{eqnarray*}
&&\vartheta\left[{ {\alpha^{(g+1)}} \atop {\beta^{(g+1)}} }\right]
\left( \int_{y}^{x} {\nu}^{(g+1)} |\Omega^{(g+1)}\right)
=\\
&&\left( -\frac{\rho}{K_0^2} \right)^{\half\kappa^2}(-\theta_{g+1})^{-\kappa}U(x,y)^{\kappa}
\vartheta\left[{ {\alpha^{(g)}} \atop {\beta^{(g)}} }\right]
\left( \int_{y}^{x} {\nu}^{(g)}+\kappa z_{p_1,p_2} |\Omega^{(g)}\right)(1+O(\rho)).
\label{eq:thetagplus}
\end{eqnarray*} 
Since  $U(x,x)=1$ we find that for $\kappa\neq-\half$
\begin{equation*}
\frac{ \vartheta \left[ {{\alpha^{(g+1)}} \atop {\beta^{(g+1)}}} \right] 
\left( \int_{y}^{x}\nu^{(g+1)}|\Omega^{(g+1)} \right) } %%
  {\vartheta \left[ {{\alpha^{(g+1)}} \atop {\beta^{(g+1)}}}\right](0|\Omega^{(g+1)})}=
  U(x,y)^{\kappa}\frac{ \vartheta \left[ {{\alpha^{(g)}} \atop {\beta^{(g)}}} \right] 
\left( \int_{y}^{x}\nu^{(g)}+\kappa z_{p_1,p_2}|\Omega^{(g)} \right) } %%
  {\vartheta \left[ {{\alpha^{(g)}} \atop {\beta^{(g)}}}\right](\kappa z_{p_1,p_2}|\Omega^{(g)})}(1+O(\rho)),
\end{equation*}
is holomorphic in $\rho$ for $|\rho|<r_{1}r_{2}$. Combining this result with \eqref{E_rhoplus} we immediately find \eqref{eq:Skappa} using the definition of the Szeg\"o kernel \eqref{Szegodefn}.

For $\kappa=-\half$ the lower bound on $|n+\alpha_{g+1}^{(g+1)}|= |\kappa|$ is satisfied for two values of $n$ so that
\begin{eqnarray*}
&&\vartheta\left[{ {\alpha^{(g+1)}} \atop {\beta^{(g+1)}} }\right]
\left( \int_{y}^{x} {\nu}^{(g+1)} |\Omega^{(g+1)}\right)
=\\
&&\left( -\frac{\rho}{K_0^2} \right)^{\frac{1}{8}}
\left[(-\theta_{g+1})^{-{\half}}U(x,y)^{\half}
\vartheta\left[{ {\alpha^{(g)}} \atop {\beta^{(g)}} }\right]
\left( \int_{y}^{x} {\nu}^{(g)} +\half z_{p_1,p_2}|\Omega^{(g)}\right)\right.\\
&&+\left.(-\theta_{g+1})^{{\half}}U(x,y)^{-\half}
\vartheta\left[{ {\alpha^{(g)}} \atop {\beta^{(g)}} }\right]
\left( \int_{y}^{x} {\nu}^{(g)}-\half z_{p_1,p_2} |\Omega^{(g)}\right)\right]\left(1+O(\rho)\right).
\end{eqnarray*} 
which eventually leads to \eqref{eq:Shalf}. \quad $\square$

\medskip
We next note that, similarly to \eqref{eq:SSpoles},   
$S_{\kappa}^{(g)}(x,z_a)S^{(g+1)}(z_a, y)$ is a meromorphic 1-form in $z_a$ periodic on the $\Sigma^{(g)}$ cycles $a_i,b_i$ for $i=1, \ldots, g$ with simple poles 
\begin{eqnarray}
S_{\kappa}^{(g)}(x,z_a)S^{(g+1)}(z_a, y)&\sim& 
\frac{dz_a}{x-z_a}S^{(g+1)}(x,y)
\mbox{ for }z_a \sim x,\notag\\
S_{\kappa}^{(g)}(x,z_a)S^{(g+1)}(z_a, y)&\sim& 
\frac{dz_a}{z_a-y}S_{\kappa}^{(g)}(x,y) 
\mbox{ for }z_a \sim y.
\label{eq:SSpoles_rho}
\end{eqnarray}   
Furthermore, $S_{\kappa}^{(g)}(x,z_a)S^{(g+1)}(z_a, y)$ is also periodic on the $a_{g+1}$ cycle defined by $\mathcal{C}_{2}(z_{2})\sim -\mathcal{C}_{1}(z_{1})$. 
This follows from applying \eqref{agplus_mult} to \eqref{eq:SgpSkappa} so that
\begin{equation}
S_{\kappa}^{(g)}(x,e^{\tpi  }z_a)=e^{\tpi \kappa(\bar{a}-a)}S_{\kappa}^{(g)}(x,z_a),
\label{eq:Skappaperiod}
\end{equation} 
(or alternatively we may apply $U(x,e^{\tpi  }z_a)^{\kappa}=e^{\tpi \kappa(\bar{a}-a)}U(x,z_a)^{\kappa}$).
Similar properties hold for $S_{\kappa}^{(g)}(x,z_a)S^{(g+1)}(z_a, y)$. 
This leads to the following analogue of Proposition \ref{Prop_inteqn1}  
%%%%%%%%%%%%%%%%%%%%%%%%%%%%%%%%%%%%%%%%%%%%%%%%%%%%%%
\begin{proposition}
\label{Prop_inteqn2}
The Szeg\"o kernel on a genus $g+1$ Riemann surface in the $\rho$-formalism  
 for $x, y \in \widehat{\Sigma}^{(g)}$ is given by
\begin{eqnarray}
S^{(g+1)}(x,y)&=&S_{\kappa}^{(g)}(x,y) 
+\sum\limits_{a=1,2} \frac{1}{\tpi }\oint\limits_{
\mathcal{C}_{a}(z_a)
} 
S_{\kappa}^{(g)}(x,z_{a} ) 
S^{(g+1)}(z_{a}, y),   \label{w2rhointeqn1}\\
&=&S_{\kappa}^{(g)}(x,y) 
-\sum\limits_{a=1,2} \frac{1}{\tpi }\oint\limits_{
\mathcal{C}_{a}(z_a)
} 
S^{(g+1)}(x,z_{a})
S_{\kappa}^{(g)}(z_{a},y ).   \label{w2rhointeqn2}
\end{eqnarray}
\end{proposition}

%%%%%%%%%%%%%%%%%%%%%%%%%%%%%%%%%%%%%%%%%%%%%%%%%%%%%%%
\noindent
{\bf Proof.} The proof follows along the same lines as Proposition \ref{Prop_inteqn1}. Let $\sigma$ be a contour on $\widehat{\Sigma}^{(g)}$ surrounding $\mathcal{A}_{a}$ and the given points $x,y\in\widehat{\Sigma}^{(g)}$ as shown in Fig.~3.   
\medskip
%%%%%%%%%%%%%%%%%%%%%%%%%%%%%%%%%%%%%%%%%%%%%%%%%%%%%%%%%%%%%%%%%%%%%% 
\begin{center}

\medskip
\begin{picture}(300,100)

\put(175,52){\qbezier(-120,20)(-80, 40)(-90,70)}%

\put(175,52){\qbezier(-90,70)(-30, 50)(-20,70)}%

\put(175,52){\qbezier(-20,70)(0, 40)(40,60)}%

\put(175,52){\qbezier(40,60)(30, 40)(70,20)}%

%%%%%%%%%%%%%%%%%%%%%%%%%%%

\put(175,52){\qbezier(-120,20)(-90, -20)(-90,-60)}%

\put(175,52){\qbezier(-90,-60)(-40, -45)(-20,-60)}%

\put(175,52){\qbezier(-20,-60)(10, -45)(40,-60)}%

\put(175,52){\qbezier(40,-60)(40, -40)(70, 20)}%

%%%%%%%%%%%%%%%%%%%%%%%%%%%%%%%%%%%%%%%%%%%%%
% left annulus 

\put(120,50){\circle{40}}
\put(120,50){\circle{22}}
\put(120,50){\makebox(0,0){$\cdot$}}

\put(120,15){\makebox(0,0){${z_{1} =0}$}}

\put(121,50){\vector(1,4){0}}%arrow

\put(185,50){\vector(1,4){0}}%arrow

{\qbezier(121,25)(115,35)(121,50)}%

 {\qbezier(185,25)(180,35)(185,50)}%

\put(130,90){\makebox(0,0){$\cdot \; x$}}

%right annulus

%%%
\put(185,50){\circle{50}}
\put(185,50){\circle{22}}
\put(206,50){\vector(1,4){0}}%arrow
\put(141,50){\vector(1,4){0}}%arrow

%%%
\put(80.5,92){\vector(-1, -1){0}}%arrow

\put(120,112.5){\vector(-1, 0){0}}%arrow

\put(180,104){\vector(-1, 0){0}}%arrow

\put(221,89){\vector(-1, 2){0}}%arrow

\put(185,50){\makebox(0,0){$\cdot$}}

%%%%%%%%%%%%%%%%%%%%%%%%%%%%%%%%%%%%%%%%%%%
\put(130,-1){\vector(1, 0){0}}%arrow

\put(190,-1){\vector(1, 0){0}}%arrow

\put(75,40){\vector(1, -2){0}}%arrow

\put(185,15){\makebox(0,0){$z_2=0$}}

% line and r2 label

\put(150,40){\makebox(0,0){$\mathcal{C}_{1} $}}

\put(215,40){\makebox(0,0){$\mathcal{C}_2$}}

\put(190,90){\makebox(0,0){$\cdot \; y$}}

\put(250,50){\makebox(0,0){$\sigma$}}

\put(235,50){\vector(1,2){0}}%arrow

%%%%%%%%%%%%%%%%%%%%%%%%%%%%%%%%%%%%%%%%%%%%%%%%%%%
% line and rho/r2 label
\put(120,50){\line(-1,2){5}}

\put(115,55){\vector(1,0){0}}%arrow

\put(10,35){\qbezier(105, 20)(85,19)(80,25)}%

\put(93,70){\makebox(0,0){$|\rho|/r_2$}}

% line and rho/r1 label
\put(185,50){\line(1,2){5}}

\put(190,57){\vector(-1, 1){0}}%arrow

\put(110,35){\qbezier(105, 25)(95, 15) (82, 20)}%

\put(215,70){\makebox(0,0){$|\rho|/r_{1} $}}

\end{picture}

\medskip
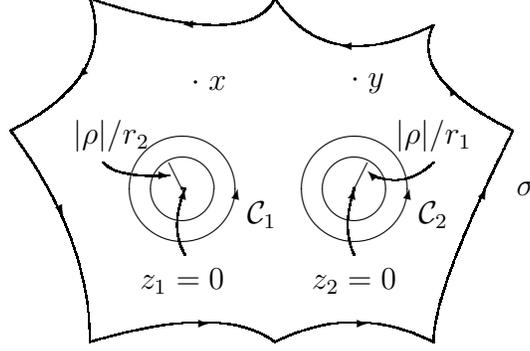
{\small Fig. 3: Contour $\sigma$}

\end{center}
Cauchy's Theorem and \eqref{eq:SSpoles_rho} imply
\begin{eqnarray*}
\label{cau0}
0&=&\frac{1} {\tpi } 
\oint_{\sigma}  
S_{\kappa}^{(g)}(x,z ) 
S^{(g+1)}(z,y)\\
&=& -S^{(g+1)}(x,y)+S_{\kappa}^{(g)}(x,y)+\sum\limits_{a=1,2} \frac{1}{\tpi }\oint\limits_{
\mathcal{C}_{a}(z_a)
} 
S_{\kappa}^{(g)}(x,z_{a} ) 
S^{(g+1)}(z_{a}, y),
\end{eqnarray*}
recalling that $S_{\kappa}^{(g)}(x,z_{a}) S^{(g+1)}(z_{a}, y)$ is periodic on $\mathcal{C}_{a}$. Thus \eqref{w2rhointeqn1} follows. A similar argument holds for \eqref{w2rhointeqn2}. \hfill $\square$ 

\medskip
%%
%%%%%%%%%%%%%%%%%%%%%%%%%%%%%%%%%%%%%%%%%%%%%%%%%%%%%

We next define weighted moments of $S^{(g+1)}(x,y)$. Let 
\begin{equation*}
\label{kakappa}
k_a =k  +(-1)^{\bar{a}}\kappa,  
\end{equation*}
for $a=1,2$ and integer $k\ge 1$ and define 
\begin{eqnarray}
&&Y_{a b}(k,l) = Y_{a b}\left[ {\theta}^{(g+1)}  \atop {\phi}^{(g+1)} \right](k,l) \nonumber\\
&&= \frac{ \rho^{\half(k_a+l_b-1) }} { (\tpi )^2 }  
\oint_{\mathcal{C}_{{\bar{a}}}(x_{\bar{a}})}
\oint_{\mathcal{C}_{b}(y_b)}
(x_{\bar{a}})^{-k_a} (y_b)^{-l_b}
S^{(g+1)} (x_{\bar{a}},y_b) 
dx_{\bar{a}}^{\half} dy_b^{\half}. 
\label{Yijdef}
\end{eqnarray} 
 We define 
$Y=\left(Y_{ab}(k,l)\right)$ 
to be the infinite matrix indexed by $a,k$ and $b,l$. 
From (\ref{Sz_skewsym}) we note the skew-symmetry property 
\begin{equation}
{Y}_{ab}\left[ {{\theta}^{(g+1)}  \atop {\phi}^{(g+1)}  }\right](k,l)= - 
{Y}_{\bar{b} \bar{a}} \left[ {(\theta^{(g+1)})^{-1}} \atop  {(\phi^{(g+1)})^{-1}}\right]
(k,l). \label{Yskew}
\end{equation}

We also introduce moments for $S^{(g)}_{\kappa}(x,y)$ 
\begin{eqnarray}
&&G_{a b}(k,l)= G_{a b}\left[ {\theta}^{(g)}  \atop {\phi}^{(g)} \right](\kappa;k,l) \nonumber\\
&&= \frac{ \rho^{\half(k_a+l_b-1) }} { (\tpi )^2 }  
\oint_{\mathcal{C}_{{\bar{a}}}(x_{\bar{a}})}
\oint_{\mathcal{C}_{b}(y_b)}
(x_{\bar{a}})^{-k_a} (y_b)^{-l_b}
S^{(g)}_{\kappa} (x_{\bar{a}},y_b) 
dx_{\bar{a}}^{\half} dy_b^{\half}, 
\label{Gijdef}
\end{eqnarray} 
with associated infinite matrix $G=\left(G_{ab}(k,l)\right)$. This also satisfies a skew-symmetry property 
\begin{equation}
{G}_{ab}\left[ {{\theta}^{(g)}  \atop {\phi}^{(g)}  }\right](\kappa;k,l)= - 
{G}_{\bar{b} \bar{a}} \left[ {(\theta^{(g)})^{-1}} \atop  {(\phi^{(g)})^{-1}}\right]
(-\kappa;k,l). \label{Gskew}
\end{equation}
Finally we define half-order differentials
\begin{eqnarray}
h_{a}(k,x)=h_{a}\left[ {\theta}^{(g)}  \atop {\phi}^{(g)} \right](\kappa;k,x) =  
\frac{\rho^{\half(k_{a} - \half)}}{\tpi }
\oint_{\mathcal{C}_{a}(y_{a})}
y_{a}^{-k_a} 
{{S}^{(g)}_{\kappa}}(x, y_{a}) 
dy_{a}^{\half},&& 
\label{hdef1}\\
\bar{h}_{a}(k,y)=\bar{h}_{a}\left[ {\theta}^{(g)}  \atop {\phi}^{(g)} \right](\kappa;k,y) =  
\frac{\rho^{\half(k_{a} - \half)}}{\tpi }
\oint_{\mathcal{C}_{\bar{a}}(x_{\bar{a}})}
x_{\bar{a}}^{-k_a} 
{{S}^{(g)}_{\kappa}}(x_{\bar{a}}, y) 
dx_{\bar{a}}^{\half},&& 
\label{hdef2}
\end{eqnarray}
and let $h(x)=(h_{a}(k, x))$ and $\bar{h}(y)=(\bar{h}_{a}(k,y))$ denote the
infinite row vectors indexed by $a,k$. These are related by skew-symmetry with
\begin{equation}
h_{a}
\left[ {\theta}^{(g)}  \atop {\phi}^{(g)} \right]
(\kappa;k,x) =-
\bar{h}_{\bar{a}} 
\left[ {(\theta^{(g)})^{-1}} \atop  {(\phi^{(g)})^{-1}}\right]
(-\kappa;k,x).
\label{hhbar2}
\end{equation}   
These moments can be inverted to obtain   
\begin{eqnarray}
S^{(g)}_{\kappa}(x,y_a) &=&
\sum_{k\geq 1} 
\rho^{-\frac{k_a}{2} +\frac{1}{4}} h_{a}(k,x)y_a^{k_a-1}dy_a^{\half}\label{hS}\\
S^{(g)}_{\kappa}(x_{\bar{a}},y) &=& \sum_{k\geq 1} 
  \rho^{-\frac{k_a}{2} +\frac{1}{4}}x_{\bar{a}}^{k_a-1}  \bar{h}_{a}(k,y) dx_{\bar{a}}^{\half}\label{hbarS}.  
\end{eqnarray}

From the sewing relation \eqref{rhosew} we have
\begin{equation}
dz_a^{\half}=(-1)^{\bar a}\xi\rho^{\half} \frac{dz_{\bar a}^{\half}}{z_{\bar a}},
\label{dz1dz2rho}
\end{equation}
for $\xi\in\{\pm \sqrt{-1}\}$. We then find in a similar way to Proposition 
\ref{Prop_hXh} that
%%
%%%%%%%%%%%%%%%%%%%%%%%%%%%%%%%%%%%%%%%%%%%%%%%%%%%%%%%%%%%%%
\begin{proposition}
\label{Prop_omg1g} For $x,y\in \widehat{\Sigma}^{(g)}$ then $S^{(g+1)}(x,y)$ is given by 
%%%%%%%%%%%%%%%%%%%%%%%%%%%%%%%%%%%%%%%%%%%%%%%%%%%%%%%%%
\begin{equation}
\label{Sgplus1}
S^{(g+1)}(x,y)= 
S_{\kappa}^{(g)}(x,y) +\xi h(x) D^{\theta}
\left( I +\xi Y D^{\theta} \right)  
\bar{h}(y)^{T}, 
\end{equation}
for infinite diagonal matrix $D^{\theta}(k,l)=
\left[
\begin{array}{cc}
\theta_{g+1}^{-1} & 0\\
0 & -\theta_{g+1}\\
\end{array}
\right]
\delta(k,l)$.
\end{proposition}
\noindent {\bf Proof.} From \eqref{w2rhointeqn1} of Proposition \ref{Prop_inteqn2} we find
%%%%%%%%%%%%%%%%%%%%%%%%%%%%%%%%%%%%%%%%%%%%
\begin{eqnarray*}
&& S^{(g+1)}(x,y)-S^{(g)}_{\kappa}(x,y) 
= \sum_{a=1,2} \frac{1}{\tpi }\oint\limits_{\mathcal{C}_{a}(z_{a})} 
S^{(g)}_{\kappa}(x,z_a) 
S^{(g+1)}(z_a, y)
\\
&=& 
\sum_{a=1,2}
\sum\limits_{ k \ge 1} 
h_a(k,x)
\frac{\rho^{-\frac{k_a}{2} + \frac{1}{4}}} {\tpi } \oint\limits_{\mathcal{C}_{a}(z_{a})} z_a^{k_a-1}  S^{(g+1)}(z_a, y) 
 dz_a^{\half}
\\
&=&
\xi \sum_{a,k}  
h_a(k,x) 
D^{\theta}_{aa}(k,k) 
\frac{\rho^{\frac{k_a}{2} - \frac{1}{4}}} {\tpi } \oint\limits_{\mathcal{C}_{\bar{a}}(z_{\bar{a}})} 
 z_{\bar{a}}^{-k_a} S^{(g+1)}(z_{\bar{a}}, y) dz_{\bar{a}}^{\half},
\end{eqnarray*}
%%%%%%%%%%%%%%%%%%%%%%%%%%%%%%%%%%%%%%%%%%%%%%%%
using respectively \eqref{hS}, \eqref{bgplus_mult} and \eqref{dz1dz2rho}. Applying \eqref{w2rhointeqn2} it follows that $S^{(g+1)}(x,y)-S^{(g)}_{\kappa}(x,y)$ is given by
\begin{eqnarray*}
&&\xi 
\sum_{a,k}
h_a(k,x) 
D_{aa}^{\theta}(k,k) 
\frac{\rho^{\frac{k_a}{2} - \frac{1}{4}}} {\tpi } \oint\limits_{\mathcal{C}_{\bar{a}}(z_{\bar{a}})} 
 z_{\bar{a}}^{-k_a} S^{(g)}_{\kappa}(z_{\bar{a}}, y) dz_{\bar{a}}^{\half}\\
 &&
 -\xi 
\sum_{a,k}
h_a(k,x) 
D_{aa}^{\theta}(k,k) 
\frac{\rho^{\frac{k_a}{2} - \frac{1}{4}}} {(\tpi )^2} \oint\limits_{\mathcal{C}_{\bar{a}}(z_{\bar{a}})} 
 z_{\bar{a}}^{-k_a} .\\
 &&\sum_{b=1,2}\;\oint\limits_{\mathcal{C}_{\bar{b}}(z_{\bar{b}})} 
 S^{(g+1)}(z_{\bar{a}}, w_{\bar{b}})S^{(g)}_{\kappa}(w_{\bar{b}}, y) 
 dz_{\bar{a}}^{\half}\\
&&=\xi h(x) D^{\theta}\bar{h}(y)^{T} 
-\xi 
\sum_{a,b,k,l}
h_a(k,x) 
D_{aa}^{\theta}(k,k). \\ 
&&\frac{\rho^{\half(k_a-l_b)}} {(\tpi )^2}
\oint\limits_{\mathcal{C}_{\bar{a}}(z_{\bar{a}})}
\oint\limits_{\mathcal{C}_{\bar{b}}(w_{\bar{b}})} 
 z_{\bar{a}}^{-k_a}w_{\bar{b}}^{l_b-1}
 S^{(g+1)}(z_{\bar{a}}, w_{\bar{b}}) 
 dz_{\bar{a}}^{\half}dw_{\bar{b}}^{\half}
 \; \bar{h}_b(l,y)\\
 &&=\xi h(x) D^{\theta}\bar{h}(y)^{T} 
-h(x) D^{\theta}YD^{\theta}\bar{h}(y)^{T},
\end{eqnarray*}
on applying \eqref{hbarS}, \eqref{bgplus_mult} and \eqref{dz1dz2rho}. Hence the result follows. \hfill $\square $

%%%%%%%%%%%%%%%%%%%%%%%%%%%%%%%%%%%%%%%%%%%%%%%%%%%%%%%%%%%%%%%%%%%%%%%%%%%%%%%%%%%%%%%%
We next compute the explicit form of ${Y}$ in terms of the weighted
moment matrix $G$ for $S_{\kappa}^{(g)}$. In particular it is convenient to define $T=\xi GD^{\theta}$. From Proposition \ref{Prop_omg1g} it follows on taking moments that
\begin{equation*}
Y=G+\xi G D^{\theta} (I+\xi Y D^{\theta})G.
\end{equation*}
This can be solved recursively to obtain $Y=\sum_{n\ge 0} T^n G$. Following a similar argument to that given for Proposition \ref{Prop_F1F2} we then find
\begin{proposition}
\label{PropYij} $Y=(I-T)^{-1}G$ where $(I-T)^{-1}=\sum_{n\ge 0} T^n$ is convergent for $|\rho |<r_{1}r_{2}$.  $\square$
\end{proposition}
This result together with Proposition \ref{Prop_omg1g} implies
\begin{theorem}
\label{Theorem_Shhrho}
$S^{(g+1)}(x,y)$ is given by 
\begin{equation*}
S^{(g+1)}(x,y)
=S^{(g)}_{\kappa}(x,y)
 +\xi h(x)D^{\theta} (I-T)^{-1}  
 \bar{h}^{T}(y).\quad \square
\end{equation*}
\end{theorem}
Finally, similarly to Theorem \ref{theorem_Det} we may define $\det \left(I-T\right)$ 
and find
\begin{theorem}
\label{theorem_Det_rho} $\det \left(I-T\right)$ is non-vanishing and holomorphic in 
$\rho$ for $|\rho |<r_{1}r_{2}$. $\square$
\end{theorem}
%%
%%%%%%%%%%%%%%%%%%%%%%%%%%%%%%%%%%%%%%%%%%%%%%%%%%%%%%%%%%%%%%%%%%%%%%%%%%%%%%%
\subsection{Self-Sewing a Sphere}
\label{self_sewing_sphere}
%%
%%%%%%%%%%%%%%%%%%%%%%%%%%%%%%%%%%%%%%%%%%%%%%%%%%%%%%%%%%%%%%%%%%%%%%%%%%%%%%%%
We consider the example of sewing the Riemann
sphere $\Sigma^{(0)}=\mathbb{C}\cup \{\infty \}$ to itself to form a torus. 
Choose local coordinates $z_{2}=z\in \mathbb{C}$ in the
neighborhood of the origin $p_2=0$, and $z_{1}=1/z^{\prime }$ for $z^{\prime }$ in
the neighborhood of the point at infinity $p_{1} =\infty$. Identify the annular regions 
$|q|r_{\bar{a}}^{-1}\leq \left\vert z_{a}\right\vert \leq r_{a}$ for a
complex sewing parameter $\rho=q$ obeying $|q|\leq r_{1}r_{2}$, via the sewing relation 
\begin{equation*}
z=qz^{\prime}.  
\label{spheresew1}
\end{equation*}
These annular regions do not intersect on the sphere provided 
$r_{1}r_{2}<1$ so that $|q|<1$.  Furthermore, the sewing relation implies $\log z=\log z^{\prime}+\tpi \tau+\tpi  k$ for integer $k$ where $q=e^{\tpi \tau}$. This is the standard parameterization of the torus with periods $\tpi \tau$ and $\tpi $ and modular parameter $\tau\in \mathbb{H}_{1} $.

We now show that the results of the previous subsection allow us to recover the genus one Szeg\"o kernel \eqref{S1} from the genus zero one. 
For $x,y\in\mathbb{C}$ the genus zero prime form and Szeg\"o kernel are given by 
\begin{eqnarray}
E^{(0)}(x, y)&=&(x-y) dx^{-\half} dy^{-\half},
\label{eq:E0}
\\
S^{(0)}(x, y)&=&\frac{1}{x-y}  dx^{\half} dy^{\half}.
\label{eq:S0}
\end{eqnarray}
Let $\theta=\theta^{(1)}$ and $\phi=\phi^{(1)}=-e^{\tpi  \kappa}$ denote the multipliers on the torus cycles.  Then since $p_{1} =\infty$ and $p_2=0$ we find $U(x,y)=x/y$ so that \eqref{eq:Skappa} and \eqref{eq:Shalf} imply
\begin{equation}
S_{\kappa}^{(0)}(x,y) =
\frac{x^{\kappa} y^{-\kappa}}{x-y}  dx^{\half} dy^{\half}
+\frac{\theta}{1-\theta}
\frac{dx^{\half} dy^{\half}}{x^{\half}y^{\half}}\delta_{\kappa,-\half}.
\label{Skappa0}
\end{equation} 
Computing moments one finds that for $\kappa\neq-\half$ the half-differentials are
\begin{eqnarray*}
h_{1}(k, x) &=& 
-\xi q^{\half(k+\kappa-\half)}  x^{k+\kappa-1 } dx^{\half},
\\
h_{2}(k, x) &=& 
q^{\half(k-\kappa-\half)}  x^{-k+\kappa}  dx^{\half},
\\
\bar{h}_{1}(k, y) &=& 
-q^{\half(k+\kappa-\half)}  y^{-k-\kappa}  dy^{\half},
\\
\bar{h}_{2}(k, y) &=& 
\xi q^{\half(k-\kappa-\half)}  y^{k-\kappa-1 } dy^{\half},
\end{eqnarray*}
for $x,y\in \widehat{\Sigma}^{(0)}$ and the moment matrix $T=\xi GD^{\theta}$ is diagonal with 
\begin{eqnarray*}
T_{ab}(k,l)=\theta^{a-\bar{a}}q^{k_a-\half}\delta_{ab}\delta(k,l).
\end{eqnarray*}

Altogether we find from Theorem \ref{Theorem_Shhrho} that for $\kappa\neq-\half$ and $x,y\in \widehat{\Sigma}^{(0)}$
\begin{eqnarray*}
S^{(1)}(x,y) 
&=&  
S_{\kappa}^{(0)}(x,y) 
+ \xi h(x)D^{\theta} (I-T)^{-1} \bar{h}^{T}(y)\\
&=&  
\left[-\frac{\left(\frac{x}{y}\right)^{\kappa + \half}}{1-\frac{x}{y}}  
- \sum\limits_{k\ge 1}  
\frac{\theta^{-1}  q^{k+\kappa - \half} } {1- \theta^{-1} q^{k+\kappa-\half} } \; \left(\frac{x}{y}\right)^{k+\kappa - \half} 
\right.\\
&&\left.+\sum\limits_{k\ge 1}  
\frac{\theta  q^{k-\kappa - \half} } {1- \theta q^{k-\kappa-\half} } \;
\left(\frac{y}{x}\right)^{k-\kappa - \half}\right]
\frac{dx^{\half}dy^{\half}}{x^{\half} y^{\half}}.
\end{eqnarray*} 
Denoting $q_u=e^u$ for any $u$ we define $X,Y$ by $x=q_X$,  $y=q_Y$ and let  $Z=X-Y$. We also define $\lambda=\kappa+\half$ with $0<\lambda<1$ and obtain  
\begin{eqnarray*}
S^{(1)}(X,Y) 
&=&   
\left[-\frac{q_Z^{\lambda}}{1-q_Z}  
- \sum\limits_{k\ge 0}  
\frac{\theta^{-1}  q^{k+\lambda} } {1- \theta^{-1} q^{k+\kappa+\half} } \; q_Z^{k+\lambda} 
\right.\\%%
&&\left.+\sum\limits_{k\le -1}  
\frac{\theta  q^{-k-\lambda} } {1- \theta q^{-k-\lambda} } \;
q_Z^{k+\lambda}\right]
dX^{\half}dY^{\half}\\
&&=-\sum\limits_{k\in\mathbb{Z}}\frac{q_Z^{k+\lambda}}{1- \theta^{-1} q^{k+\lambda}}dX^{\half}dY^{\half}=P_{1} (Z,q)dX^{\half}dY^{\half},
\end{eqnarray*}
from \eqref{P1}. A similar result also holds for $\kappa=-\half$ for $\theta\neq 0$ i.e. $(\theta,\phi)\neq (0,0)$.

Lastly, we note that $(I-T)^{-1}$ is convergent for $|q|<1$ and that furthermore
\begin{equation}
\label{detTg1}
\det(I-T)=\prod_{k \ge 1} 
\left( 1-\theta^{-1} q^{k + \kappa - \half }\right)\; 
\left(1-\theta\ q^{k - \kappa-\half}\right), 
\end{equation}
is holomorphic for $|q|<1$ from Theorem  
\ref{theorem_Det_rho}.
In vertex operator algebra theory, $\det(I-T)$ is related to the genus one partition function for a continuous orbfolding of a rank two free fermion system e.g. \cite{MTZ}. Furthermore, the infinite product \eqref{detTg1} is part of that arising in the Jacobi triple identity on applying the bosonic decomposition of this theory. 

\medskip
%%%%%%%%%%%%%%%%%%%%%%%%%%%%%%%%%%%%%%%%%%%%%%%%%%%%%%%%%%%%%%%%%%%%%%%%
%%%%%%%%%%%%%%%%%%%%%%%%%%%%%%%%%%%%%%%%%%%%%%%%%%%%%%%%%%%%%%%%%%%%%%%%
\subsection{Self-Sewing a Torus}
\label{Sect_Rho_Torus}
We next consider the example of self-sewing an oriented torus
$\Sigma^{(1)}=\mathbb{C}/{\Lambda} $ for lattice ${\Lambda}=\tpi (\mathbb{Z}\tau\oplus \mathbb{Z})$ and $\tau\in \mathbb{H}_{1}$. This is discussed in detail in ref. \cite{MT1}.
Define annuli $\mathcal{A}_{a},a=1,2$ centered at $p_{1} =0$ and $p_2=w$ of  $\Sigma^{(1)}$ with local
coordinates $z_{1}=z$ and $z_{2}=z-w$ respectively. Take the outer radius of 
$\mathcal{A}_{a}$ to be $r_{a}<\half D(q)$ for $D(q)=\min_{\lambda \in {\Lambda}, \lambda \neq 0}|\lambda |$ and the inner radius
to be $|{\rho }|/r_{\bar{a}}$, with $|\rho|\leq r_{1}r_{2}$.   
Identifying the annuli via \eqref{rhosew} we obtain 
a compact genus two Riemann surface $\Sigma^{(2)}$
parameterized by 
\begin{equation}
\mathcal{D}^{\rho }=\{(\tau ,w,\rho )\in \mathbb{H}_{1}\times \mathbb{C}
\times \mathbb{C}\ :\ |w-\lambda |>2|\rho |^{\half}>0,\ \lambda \in {\Lambda}\}.
  \label{Drho}
\end{equation}
%%%%%%%%%%%%%%%%%%%%%%%%%%%%%%%%%%%%%%%%%%%%%%%%%%%%%%%%%%%%%%%%
%%
For $x,y\in\Sigma^{(1)}$ the genus one prime form and Szeg\"o kernel with multipliers $\theta_{1} =-e^{-\tpi \beta_{1}}$ and $\phi_{1} =-e^{\tpi \alpha_{1}}$ are given by \eqref{Ktheta} and \eqref{S1}. 
Let $\theta_{2}=-e^{-\tpi \beta_{2}}$ and $\phi_{2}=-e^{\tpi \alpha_{2}}=-e^{\tpi  \kappa}$ denote the multipliers on $a_{2},b_{2}$ cycles. Then, in this case
\begin{equation*}
U(x,y)=\frac {\vartheta_{1} (x-w,\tau) \vartheta_{1} (y,\tau)} 
{ \vartheta_{1} (x,\tau)  \vartheta_{1} (y-w,\tau)  },
\end{equation*}
and $z_{0,w}=\kappa w$
so that for $\kappa\neq -\half$
\begin{equation*}
S_{\kappa}^{(1)}
\left[ {\theta_{1} } \atop {\phi_{1} }\right]
(x,y\vert\tau,w)=
\left(
\frac {\vartheta_{1} (x-w,\tau) \vartheta_{1} (y,\tau)} 
{ \vartheta_{1} (x,\tau)  \vartheta_{1} (y-w,\tau)}
\right)^{\kappa}
\frac{ \vartheta \left[ {\alpha_{1} } \atop {\beta_{1} } \right] 
\left( x-y +\kappa w,\tau\right)}
{\vartheta \left[ {\alpha_{1} } \atop {\beta_{1} }\right] 
\left(\kappa w,\tau \right)  K(x-y,\tau)}
dx^{\half} dy^{\half},
\end{equation*}
with a similar result for $\kappa= -\half$. We take $\kappa\neq-\half$ from now on.

It is straightforward to see that 
\begin{equation}
S_{\kappa}^{(1)}\left[ {\theta_{1} } \atop {\phi_{1} }\right](x,y\vert\tau,w) =
  S_{-\kappa}^{(1)}\left[ {\theta_{1} } \atop {\phi_{1} }\right](x-w,y-w\vert\tau,-w).
 \label{Sksym}
\end{equation} 
Computing moments and using \eqref{hhbar2} and \eqref{Sksym} the half-differentials \eqref{hdef1}, \eqref{hdef2} 
for $x\in \widehat{\Sigma}^{(1)}$ and $\kappa\neq-\half$ are given by 
\begin{eqnarray}
h_{1}
\left[ {\theta_{1} } \atop {\phi_{1} }\right]
(\kappa;k, x\vert\tau,w,\rho) &=& 
\frac{\rho^{\half (k+\kappa -\half)}}{\tpi}
\left(
\frac {\vartheta_{1} (x-w,\tau)}{\vartheta_{1} (x,\tau)}
\right)^{\kappa}
\frac{dx^\half}{\vartheta \left[ {\alpha_{1} } \atop {\beta_{1} }\right] 
\left(\kappa w,\tau \right)} \notag\\
&&
\oint_{{\cal C}_{1} (y)}
y^{-k-\kappa}
\left(
\frac {\vartheta_{1} (y,\tau)} {\vartheta_{1} (y-w,\tau)}
\right)^{\kappa}
\frac{ \vartheta \left[ {\alpha_{1} } \atop {\beta_{1} } \right] 
\left( x-y +\kappa w,\tau\right)}
{K(x-y,\tau)}dy, \notag\\
h_{2}
\left[ {\theta_{1} } \atop {\phi_{1} }\right]
(\kappa;k, x\vert\tau,w,\rho) &=&
h_{1}
\left[ {\theta_{1} } \atop {\phi_{1} }\right]
(-\kappa;k, x-w\vert\tau,-w,\rho),\notag \\
\bar{h}_{1}
\left[ {\theta_{1} } \atop {\phi_{1} }\right]
(\kappa;k, x\vert\tau,w,\rho) &=&
-h_{1}
\left[ {\theta_{1} ^{-1}} \atop {\phi_{1} ^{-1}}\right]
(\kappa;k, x-w\vert\tau,-w,\rho), \notag \\
\bar{h}_{2}
\left[ {\theta_{1} } \atop {\phi_{1} }\right]
(\kappa;k, x\vert\tau,w,\rho) &=&
-h_{1}
\left[ {\theta_{1} ^{-1}} \atop {\phi_{1} ^{-1}}\right]
(-\kappa;k, x\vert\tau,w,\rho). \label{htorus}
\end{eqnarray}
Similarly, using \eqref{Sksym}, the moment matrix \eqref{Gijdef} is given by
\begin{eqnarray}
G_{11}
\left[ {\theta_{1} } \atop {\phi_{1} }\right]
(\kappa;k, l\vert\tau,w,\rho) &=& 
\frac{\rho^{\kappa+\half (k+l-1)}}{(\tpi)^2}
\oint_{{\cal C}_{2} (x_2)}
\oint_{{\cal C}_{1} (y_1)}{x_2}^{-k-\kappa} {y_1}^{-l-\kappa}\notag\\
&&
S_{\kappa}^{(1)}\left[ {\theta_{1} } \atop {\phi_{1} }\right](x_2,y_1\vert\tau,w)
dx_{2}^\half dy_{1}^\half,
\notag\\
&&=G_{22}
\left[ {\theta_{1} } \atop {\phi_{1} }\right]
(-\kappa;k, l\vert\tau,-w,\rho),\notag\\
G_{21}
\left[ {\theta_{1} } \atop {\phi_{1} }\right]
(\kappa;k, l\vert\tau,w,\rho) &=& 
\frac{\rho^{\half (k+l-1)}}{(\tpi)^2}
\oint_{{\cal C}_{1} (x_1)}
\oint_{{\cal C}_{1} (y_1)}{x_1}^{-k+\kappa} {y_1}^{-l-\kappa}\notag\\
&&
S_{\kappa}^{(1)}\left[ {\theta_{1} } \atop {\phi_{1} }\right](x_1,y_1\vert\tau,w)
dx_{1}^\half dy_{1}^\half,\notag\\
&&=G_{12}
\left[ {\theta_{1} } \atop {\phi_{1} }\right]
(-\kappa;k, l\vert\tau,-w,\rho).\label{Gtorus}
\end{eqnarray}
The genus two Szego kernel is determined
for $T=\xi G\left[ {\theta_{1} } \atop {\phi_{1} }\right] D^{\theta_2}$ by \eqref{Sgplus1}
\begin{eqnarray}
\label{Sg2}
&S^{(2)}
\left[
{\theta_1} \atop {\theta_2}	
\right]
(x,y\vert\tau,w,\rho)
=&\notag \\
&S_{\kappa}^{(1)}\left[ {\theta_{1} } \atop {\phi_{1} }\right](x,y\vert\tau,w) +
 \xi h\left[ {\theta_{1} } \atop {\phi_{1} }\right](x) D^{\theta_2}
\left( I - T \right)^{-1}  
\bar{h}^{T}\left[ {\theta_{1} } \atop {\phi_{1} }\right](y).& 
\end{eqnarray} 

\subsubsection{Modular Invariance}
We now consider the modular invariance of \eqref{Sg2} under the action of a particular subgroup $L\subset Sp(4,{\mathbb Z})$ and verify that \eqref{Szmod} holds.
We define $L$ as follows \cite{MT1}. Consider  $\hat{H} \subset Sp(4,{\mathbb Z})$ with elements
\begin{equation}
\mu (a,b,c)=\left( 
\begin{array}{cccc}
1 & 0 & 0 & b \\ 
a & 1 & b & c \\ 
0 & 0 & 1 & -a \\ 
0 & 0 & 0 & 1%
\end{array}%
\right).  \label{mudef}
\end{equation}%
$\hat{H}$ is generated by $A=\mu (1,0,0)$, $B=\mu (0,1,0)$ and $C=\mu (0,0,1)$ with relations $[A,B]C^{-2}=[A,C]=[B,C]=1$. We also define $\Gamma _{1}\subset Sp(4,\mathbb{Z})$ where $\Gamma_{1}\cong SL(2,\mathbb{Z})$ with elements
\begin{equation}
\gamma_1 =\left( 
\begin{array}{cccc}
a_{1} & 0 & b_{1} & 0 \\ 
0 & 1 & 0 & 0 \\ 
c_{1} & 0 & d_{1} & 0 \\ 
0 & 0 & 0 & 1%
\end{array}%
\right),\ \ a_{1}d_{1}-b_{1}c_{1}=1.
\label{gamma}
\end{equation}%
Together these groups generate $L=\hat{H}\rtimes\Gamma _{1}\subset Sp(4,\mathbb{Z})$ with center $Z(L)=\langle C\rangle $ where $J=L/Z(L)\cong 
\mathbb{Z}^{2}\rtimes SL(2,\mathbb{Z})$ is the Jacobi group. 

From Lemma~15 of \cite{MT1} we find that $L$ acts on the domain $\mathcal{D}^{\rho }$ of \eqref{Drho} as follows: 
\begin{eqnarray}
\mu(a, b, c).(\tau, w, \rho) &=& (\tau, w+2 \pi i a \tau + 2 \pi i b, \rho),
\label{muaction} \\
\gamma_1. (\tau, w, \rho) &=& \left( \frac{a_{1} \tau+b_{1}}{c_{1}\tau + d_{1}}, \frac{w}{c_{1}
\tau + d_{1}}, \frac{\rho}{(c_{1} \tau + d_{1})^2} \right).  \label{gammaonDrhoaction}
\end{eqnarray}
The kernel of the action is $Z(L)$, so that the effective action is that of 
$J$. However,  this action is lifted to $L$ when considering the covering space $\widehat{\mathcal{D}}^{\rho}$  of $\mathcal{D}^{\rho }$ 
for which $\Omega^{(g+1)}_{g+1,g+1}$ of \eqref{Omgpgp_rho} is single-valued (Theorems 10, 11 of \cite{MT1}). In particular, one finds that $C$ acts as 
\begin{eqnarray}
C.(\tau, w, \rho) =(\tau, w, e^{\tpi}\rho),
\label{Caction}
\end{eqnarray}
which has a non-trivial action on $\widehat{\mathcal{D}}^{\rho }$.

Let us now consider the action of $L$ on $S^{(2)}
\left[
{\theta^{(2)}} \atop {\phi^{(2)}}	
\right]
(x,y\vert\tau,w,\rho)$. This is partly determined by the action of $J$ on $S_{\kappa}^{(1)}\left[ {\theta_{1} } \atop {\phi_{1} }\right](x,y\vert w,\tau)$. For $\gamma_1\in \Gamma_1$ it is clear from \eqref{eq:S1mod} that
\begin{equation*}
S_{\kappa}^{(1)}
\left[ {{\theta_{1}}^a \phi_{1}^b } \atop {{\theta_{1}}^c \phi_{1}^d}\right]
(\gamma_1 x,\gamma_1 y\vert\gamma_1(\tau,w))=S_{\kappa}^{(1)}
\left[ {\theta_{1} } \atop {\phi_{1} }\right]
(x,y\vert\tau,w).
\end{equation*} 
$h_{a}\left[ {\theta_{1} } \atop {\phi_{1} }\right]
(\kappa;k, x\vert\tau,w,\rho)$ and $G_{ab}\left[ {\theta_{1} } \atop {\phi_{1} }\right]
(\kappa;k, l\vert\tau,w,\rho)$ are similarly $\Gamma_1$ invariant so that $S^{(2)}
\left[
{\theta^{(2)}} \atop {\phi^{(2)}}	
\right]
(x,y\vert\tau,w,\rho)$ is $\Gamma_1$ invariant in a similar fashion to \eqref{eq:S2mod}.

Next we consider the action of the generators $A$, $B$ and $C$. We firstly note that \eqref{albetatilde} implies
\begin{equation}
A 
\left[
{{\theta_1} \atop {\theta_{2}}}
\atop
{{\phi_1} \atop {\phi_{2}}}	
\right]=
\left[
{{\theta_1} \atop {-\theta_2 \theta_1}}
\atop
{{-\phi_{1} \phi_{2}^{-1}} \atop {\phi_{2}}}	
\right]
,\quad
B
\left[
{{\theta_1} \atop {\theta_{2}}}
\atop
{{\phi_1} \atop {\phi_{2}}}	
\right]=
\left[
{{-\theta_{1}\phi_2} \atop {-\theta_2\phi_1}}
\atop
{{\phi_1} \atop {\phi_{2}}}	
\right]
,\quad
C
\left[
{{\theta_1} \atop {\theta_{2}}}
\atop
{{\phi_1} \atop {\phi_{2}}}	
\right]=
\left[
{{\theta_1} \atop {-\theta_{2}\phi_2}}
\atop
{{\phi_1} \atop {\phi_{2}}}	
\right].
\label{ABC}
\end{equation}  
Using  \eqref{theta} and recalling that $\phi_2=-e^{\tpi \kappa}$ we find 
\begin{eqnarray}
S_{\kappa}^{(1)}
\left[ {\theta_{1} } \atop {\phi_{1} }\right]
(x,y\vert\tau,w)
&=&
S_{\kappa}^{(1)}
\left[ {\theta_{1} } \atop {-\phi_{1} \phi^{-1}_{2}}\right]
(x,y\vert\tau,w+\tpi \tau)\label{SA}\\
&=&
S_{\kappa}^{(1)}
\left[ -{\theta_{1}\phi_{2} } \atop {\phi_{1} }\right]
(x,y\vert\tau,w+\tpi),\label{SB}
\end{eqnarray}
where the multipliers comply with those of \eqref{ABC} for $A$ and $B$ respectively.
Define infinite diagonal matrices 
\begin{equation}
E^{\alpha}(k,l)=
\left[
\begin{array}{cc}
1 & 0\\
0 & -\alpha\\
\end{array}
\right]
\delta(k,l), \quad
F^{\alpha}(k,l)=
\left[
\begin{array}{cc}
-\alpha^{-1} & 0\\
0 & 1\\
\end{array}
\right]
\delta(k,l),
\label{EFtheta}
\end{equation}
for $\alpha\in U(1)$. 
Then \eqref{htorus}, \eqref{SA} and \eqref{SB} imply
\begin{eqnarray*}
h\left[ {\theta_{1} } \atop {\phi_{1} }\right]
(\kappa;x\vert\tau,w,\rho)
&=&
h\left[ {\theta_{1} } \atop {-\phi_{1} \phi_{2}^{-1}}\right]
(\kappa;x\vert\tau,w+\tpi \tau,\rho)E^{\theta_1} \\
&=&
h\left[ -{\theta_{1}\phi_{2} } \atop {\phi_{1} }\right]
(\kappa;x\vert\tau,w+\tpi,\rho)E^{\phi_1}\\
&=&
e^{-i\pi \kappa}
h\left[ {\theta_{1}} \atop {\phi_{1} }\right]
(\kappa;x\vert\tau,w,e^{\tpi}\rho)E^{\phi_2},\\
\bar{h}\left[ {\theta_{1} } \atop {\phi_{1} }\right]
(\kappa;x\vert\tau,w,\rho)
&=&
h\left[ {\theta_{1} } \atop {-\phi_{1} \phi^{-1}_{2}}\right]
(\kappa;x\vert\tau,w+\tpi \tau,\rho)F^{\theta_1} \\
&=&
h\left[ -{\theta_{1}\phi_{2} } \atop {\phi_{1} }\right]
(\kappa;x\vert\tau,w+\tpi,\rho)F^{\phi_1}\\
&=&
e^{i\pi \kappa}
h\left[ {\theta_{1}} \atop {\phi_{1} }\right]
(\kappa;x\vert\tau,w,e^{\tpi}\rho)F^{\phi_2}.
\end{eqnarray*}
Similarly, from  \eqref{Gtorus} we find
\begin{eqnarray*}
G\left[ {\theta_{1} } \atop {\phi_{1} }\right]
(\kappa \vert\tau,w,\rho)
&=&
F^{\theta_1}G\left[ {\theta_{1} } \atop {-\phi_{1} \phi^{-1}_{2}}\right]
(\kappa \vert\tau,w+\tpi \tau,\rho)E^{\theta_1}\\
&=&
F^{\phi_1}G\left[ -{\theta_{1}\phi_{2} } \atop {\phi_{1} }\right]
(\kappa \vert\tau,w+\tpi,\rho)E^{\phi_1}\\
&=&
F^{\phi_2}G\left[ -{\theta_{1} } \atop {\phi_{1} }\right]
(\kappa \vert\tau,w,e^{\tpi}\rho)E^{\phi_2}.
\end{eqnarray*}
Noting that $E^{\alpha}D^{\theta_2}F^{\alpha}=D^{-\alpha\theta_2}$ for $\alpha=\theta_1, \phi_1$ and $\phi_2$ we may then easily confirm that $S^{(2)}
\left[
{\theta^{(2)}} \atop {\phi^{(2)}}	
\right]
(x,y\vert\tau,w,\rho)$ is invariant under the generators $A$, $B$ and $C$ respectively. Therefore $S^{(2)}
\left[
{\theta^{(2)}} \atop {\phi^{(2)}}	
\right]
(x,y\vert\tau,w,\rho)$ is modular invariant under $L$. Furthermore, since $\det(E^{\alpha}F^{\alpha})=1$ it follows that $\det \left(I-T\right)$ is also $L$ invariant.

%%%%%%%%%%%%%%%%%%%%%%%%%%%%%%%%%%%%%%%%%%%%%%%%%%%%%%%%%%%%%%%%%%

\end{document}